\documentclass{amsart}
\numberwithin{equation}{section}

\usepackage{graphicx,amsthm}
\usepackage{verbatim,amsmath,amscd,amssymb,color}
\usepackage[mathscr]{eucal}
\input xy
\xyoption{all}
\begin{document}
\renewcommand{\labelenumi}{$($\roman{enumi}$)$}
\renewcommand{\labelenumii}{$(${\rm \alph{enumii}}$)$}
\font\germ=eufm10
\newcommand{\cI}{{\mathcal I}}
\newcommand{\cA}{{\mathcal A}}
\newcommand{\cB}{{\mathcal B}}
\newcommand{\cC}{{\mathcal C}}
\newcommand{\cD}{{\mathcal D}}
\newcommand{\cE}{{\mathcal E}}
\newcommand{\cF}{{\mathcal F}}
\newcommand{\cG}{{\mathcal G}}
\newcommand{\cH}{{\mathcal H}}
\newcommand{\cK}{{\mathcal K}}
\newcommand{\cL}{{\mathcal L}}
\newcommand{\cM}{{\mathcal M}}
\newcommand{\cN}{{\mathcal N}}
\newcommand{\cO}{{\mathcal O}}
\newcommand{\cR}{{\mathcal R}}
\newcommand{\cS}{{\mathcal S}}
\newcommand{\cT}{{\mathcal T}}
\newcommand{\cV}{{\mathcal V}}
\newcommand{\cX}{{\mathcal X}}
\newcommand{\cY}{{\mathcal Y}}
\newcommand{\fra}{\mathfrak a}
\newcommand{\frb}{\mathfrak b}
\newcommand{\frc}{\mathfrak c}
\newcommand{\frd}{\mathfrak d}
\newcommand{\fre}{\mathfrak e}
\newcommand{\frf}{\mathfrak f}
\newcommand{\frg}{\mathfrak g}
\newcommand{\frh}{\mathfrak h}
\newcommand{\fri}{\mathfrak i}
\newcommand{\frj}{\mathfrak j}
\newcommand{\frk}{\mathfrak k}
\newcommand{\frI}{\mathfrak I}
\newcommand{\fm}{\mathfrak m}
\newcommand{\frn}{\mathfrak n}
\newcommand{\frp}{\mathfrak p}
\newcommand{\fq}{\mathfrak q}
\newcommand{\frr}{\mathfrak r}
\newcommand{\frs}{\mathfrak s}
\newcommand{\frt}{\mathfrak t}
\newcommand{\fru}{\mathfrak u}
\newcommand{\frA}{\mathfrak A}
\newcommand{\frB}{\mathfrak B}
\newcommand{\frF}{\mathfrak F}
\newcommand{\frG}{\mathfrak G}
\newcommand{\frH}{\mathfrak H}
\newcommand{\frJ}{\mathfrak J}
\newcommand{\frN}{\mathfrak N}
\newcommand{\frP}{\mathfrak P}
\newcommand{\frT}{\mathfrak T}
\newcommand{\frU}{\mathfrak U}
\newcommand{\frV}{\mathfrak V}
\newcommand{\frX}{\mathfrak X}
\newcommand{\frY}{\mathfrak Y}
\newcommand{\fry}{\mathfrak y}
\newcommand{\frZ}{\mathfrak Z}
\newcommand{\frx}{\mathfrak x}
\newcommand{\rA}{\mathrm{A}}
\newcommand{\rC}{\mathrm{C}}
\newcommand{\rd}{\mathrm{d}}
\newcommand{\rB}{\mathrm{B}}
\newcommand{\rD}{\mathrm{D}}
\newcommand{\rE}{\mathrm{E}}
\newcommand{\rH}{\mathrm{H}}
\newcommand{\rK}{\mathrm{K}}
\newcommand{\rL}{\mathrm{L}}
\newcommand{\rM}{\mathrm{M}}
\newcommand{\rN}{\mathrm{N}}
\newcommand{\rR}{\mathrm{R}}
\newcommand{\rT}{\mathrm{T}}
\newcommand{\rZ}{\mathrm{Z}}
\newcommand{\bbA}{\mathbb A}
\newcommand{\bbB}{\mathbb B}
\newcommand{\bbC}{\mathbb C}
\newcommand{\bbG}{\mathbb G}
\newcommand{\bbF}{\mathbb F}
\newcommand{\bbH}{\mathbb H}
\newcommand{\bbP}{\mathbb P}
\newcommand{\bbN}{\mathbb N}
\newcommand{\bbQ}{\mathbb Q}
\newcommand{\bbR}{\mathbb R}
\newcommand{\bbV}{\mathbb V}
\newcommand{\bbZ}{\mathbb Z}
\newcommand{\adj}{\operatorname{adj}}
\newcommand{\Ad}{\mathrm{Ad}}
\newcommand{\Ann}{\mathrm{Ann}}
\newcommand{\rcris}{\mathrm{cris}}
\newcommand{\ch}{\mathrm{ch}}
\newcommand{\coker}{\mathrm{coker}}
\newcommand{\diag}{\mathrm{diag}}
\newcommand{\Diff}{\mathrm{Diff}}
\newcommand{\Dist}{\mathrm{Dist}}
\newcommand{\rDR}{\mathrm{DR}}
\newcommand{\ev}{\mathrm{ev}}
\newcommand{\Ext}{\mathrm{Ext}}
\newcommand{\cExt}{\mathcal{E}xt}
\newcommand{\fin}{\mathrm{fin}}
\newcommand{\Frac}{\mathrm{Frac}}
\newcommand{\GL}{\mathrm{GL}}
\newcommand{\Hom}{\mathrm{Hom}}
\newcommand{\hd}{\mathrm{hd}}
\newcommand{\rht}{\mathrm{ht}}
\newcommand{\id}{\mathrm{id}}
\newcommand{\im}{\mathrm{im}}
\newcommand{\inc}{\mathrm{inc}}
\newcommand{\ind}{\mathrm{ind}}
\newcommand{\coind}{\mathrm{coind}}
\newcommand{\Lie}{\mathrm{Lie}}
\newcommand{\Max}{\mathrm{Max}}
\newcommand{\mult}{\mathrm{mult}}
\newcommand{\op}{\mathrm{op}}
\newcommand{\ord}{\mathrm{ord}}
\newcommand{\pt}{\mathrm{pt}}
\newcommand{\qt}{\mathrm{qt}}
\newcommand{\rad}{\mathrm{rad}}
\newcommand{\res}{\mathrm{res}}
\newcommand{\rgt}{\mathrm{rgt}}
\newcommand{\rk}{\mathrm{rk}}
\newcommand{\SL}{\mathrm{SL}}
\newcommand{\soc}{\mathrm{soc}}
\newcommand{\Spec}{\mathrm{Spec}}
\newcommand{\St}{\mathrm{St}}
\newcommand{\supp}{\mathrm{supp}}
\newcommand{\Tor}{\mathrm{Tor}}
\newcommand{\Tr}{\mathrm{Tr}}
\newcommand{\wt}{\mathrm{wt}}
\newcommand{\Ab}{\mathbf{Ab}}
\newcommand{\Alg}{\mathbf{Alg}}
\newcommand{\Grp}{\mathbf{Grp}}
\newcommand{\Mod}{\mathbf{Mod}}
\newcommand{\Sch}{\mathbf{Sch}}\newcommand{\bfmod}{{\bf mod}}
\newcommand{\Qc}{\mathbf{Qc}}
\newcommand{\Rng}{\mathbf{Rng}}
\newcommand{\Top}{\mathbf{Top}}
\newcommand{\Var}{\mathbf{Var}}
\newcommand{\pmbx}{\pmb x}
\newcommand{\pmby}{\pmb y}
\newcommand{\gromega}{\langle\omega\rangle}
\newcommand{\lbr}{\begin{bmatrix}}
\newcommand{\rbr}{\end{bmatrix}}
\newcommand{\cd}{commutative diagram }
\newcommand{\SpS}{spectral sequence}
\newcommand\C{\mathbb C}
\newcommand\hh{{\hat{H}}}
\newcommand\eh{{\hat{E}}}
\newcommand\F{\mathbb F}
\newcommand\fh{{\hat{F}}}
\newcommand\Z{{\mathbb Z}}
\newcommand\Zn{\Z_{\geq0}}
\newcommand\et[1]{\tilde{e}_{#1}}
\newcommand\ft[1]{\tilde{f}_{#1}}

\def\ge{\frg}
\def\AA{{\mathcal A}}
\def\al{\alpha}
\def\bq{B_q(\ge)}
\def\bqm{B_q^-(\ge)}
\def\bqz{B_q^0(\ge)}
\def\bqp{B_q^+(\ge)}
\def\beneme{\begin{enumerate}}
\def\beq{\begin{equation}}
\def\beqn{\begin{eqnarray}}
\def\beqnn{\begin{eqnarray*}}
\def\bfi{{\mathbf i}}
\def\bfii0{{\bf i_0}}
\def\bigsl{{\hbox{\fontD \char'54}}}
\def\bbra#1,#2,#3{\left\{\begin{array}{c}\hspace{-5pt}
#1;#2\\ \hspace{-5pt}#3\end{array}\hspace{-5pt}\right\}}
\def\cd{\cdots}
\def\ch(#1,#2){c_{#2,#1}^{-h_{#1}}}
\def\cc(#1,#2){c_{#2,#1}}
\def\CC{\mathbb{C}}
\def\CBL{\cB_L(\TY(B,1,n+1))}
\def\CBM{\cB_M(\TY(B,1,n+1))}
\def\CVL{\cV_L(\TY(D,1,n+1))}
\def\CVM{\cV_M(\TY(D,1,n+1))}
\def\ddd{\hbox{\germ D}}
\def\del{\delta}
\def\Del{\Delta}
\def\Delr{\Delta^{(r)}}
\def\Dell{\Delta^{(l)}}
\def\Delb{\Delta^{(b)}}
\def\Deli{\Delta^{(i)}}
\def\Delre{\Delta^{\rm re}}
\def\ei{e_i}
\def\eit{\tilde{e}_i}
\def\eneme{\end{enumerate}}
\def\ep{\epsilon}
\def\eeq{\end{equation}}
\def\eeqn{\end{eqnarray}}
\def\eeqnn{\end{eqnarray*}}
\def\fit{\tilde{f}_i}
\def\FF{{\rm F}}
\def\ft{\tilde{f}_}
\def\gau#1,#2{\left[\begin{array}{c}\hspace{-5pt}#1\\
\hspace{-5pt}#2\end{array}\hspace{-5pt}\right]}
\def\gl{\hbox{\germ gl}}
\def\hom{{\hbox{Hom}}}
\def\ify{\infty}
\def\io{\iota}
\def\kp{k^{(+)}}
\def\km{k^{(-)}}
\def\llra{\relbar\joinrel\relbar\joinrel\relbar\joinrel\rightarrow}
\def\lan{\langle}
\def\lar{\longrightarrow}
\def\lm{\lambda}
\def\Lm{\Lambda}
\def\mapright#1{\smash{\mathop{\longrightarrow}\limits^{#1}}}
\def\Mapright#1{\smash{\mathop{\Longrightarrow}\limits^{#1}}}
\def\mm{{\bf{\rm m}}}
\def\nd{\noindent}
\def\nn{\nonumber}
\def\nnn{\hbox{\germ n}}
\def\catob{{\mathcal O}(B)}
\def\oint{{\mathcal O}_{\rm int}(\ge)}
\def\ot{\otimes}
\def\op{\oplus}
\def\opi{\ovl\pi_{\lm}}
\def\osigma{\ovl\sigma}
\def\ovl{\overline}
\def\plm{\Psi^{(\lm)}_{\io}}
\def\qq{\qquad}
\def\q{\quad}
\def\qed{\hfill\framebox[2mm]{}}
\def\QQ{\mathbb Q}
\def\qi{q_i}
\def\qii{q_i^{-1}}
\def\ra{\rightarrow}
\def\ran{\rangle}
\def\rlm{r_{\lm}}
\def\ssl{\hbox{\germ sl}}
\def\slh{\widehat{\ssl_2}}
\def\ti{t_i}
\def\tii{t_i^{-1}}
\def\til{\tilde}
\def\tm{\times}
\def\tt{\frt}
\def\TY(#1,#2,#3){#1^{(#2)}_{#3}}
\def\ua{U_{\AA}}
\def\ue{U_{\vep}}
\def\uq{U_q(\ge)}
\def\uqp{U'_q(\ge)}
\def\ufin{U^{\rm fin}_{\vep}}
\def\ufinp{(U^{\rm fin}_{\vep})^+}
\def\ufinm{(U^{\rm fin}_{\vep})^-}
\def\ufinz{(U^{\rm fin}_{\vep})^0}
\def\uqm{U^-_q(\ge)}
\def\uqmq{{U^-_q(\ge)}_{\bf Q}}
\def\uqpm{U^{\pm}_q(\ge)}
\def\uqq{U_{\bf Q}^-(\ge)}
\def\uqz{U^-_{\bf Z}(\ge)}
\def\ures{U^{\rm res}_{\AA}}
\def\urese{U^{\rm res}_{\vep}}
\def\uresez{U^{\rm res}_{\vep,\ZZ}}
\def\util{\widetilde\uq}
\def\uup{U^{\geq}}
\def\ulow{U^{\leq}}
\def\bup{B^{\geq}}
\def\blow{\ovl B^{\leq}}
\def\vep{\varepsilon}
\def\vp{\varphi}
\def\vpi{\varphi^{-1}}
\def\VV{{\mathcal V}}
\def\xii{\xi^{(i)}}
\def\Xiioi{\Xi_{\io}^{(i)}}
\def\xxi(#1,#2,#3){\displaystyle {}^{#1}\Xi^{(#2)}_{#3}}
\def\xx(#1,#2){\displaystyle {}^{#1}\Xi_{#2}}
\def\W1{W(\varpi_1)}
\def\WW{{\mathcal W}}
\def\wt{{\rm wt}}
\def\wtil{\widetilde}
\def\what{\widehat}
\def\wpi{\widehat\pi_{\lm}}
\def\ZZ{\mathbb Z}
\def\RR{\mathbb R}

\def\m@th{\mathsurround=0pt}
\def\fsquare(#1,#2){
\hbox{\vrule$\hskip-0.4pt\vcenter to #1{\normalbaselines\m@th
\hrule\vfil\hbox to #1{\hfill$\scriptstyle #2$\hfill}\vfil\hrule}$\hskip-0.4pt
\vrule}}

\newtheorem{pro}{Proposition}[section]
\newtheorem{thm}[pro]{Theorem}
\newtheorem{lem}[pro]{Lemma}
\newtheorem{ex}[pro]{Example}
\newtheorem{cor}[pro]{Corollary}
\newtheorem{conj}[pro]{Conjecture}
\theoremstyle{definition}
\newtheorem{df}[pro]{Definition}

\newcommand{\cmt}{\marginpar}
\newcommand{\seteq}{\mathbin{:=}}
\newcommand{\cl}{\colon}
\newcommand{\be}{\begin{enumerate}}
\newcommand{\ee}{\end{enumerate}}
\newcommand{\bnum}{\be[{\rm (i)}]}
\newcommand{\enum}{\ee}
\newcommand{\ro}{{\rm(}}
\newcommand{\rf}{{\rm)}}
\newcommand{\set}[2]{\left\{#1\,\vert\,#2\right\}}
\newcommand{\sbigoplus}{{\mbox{\small{$\bigoplus$}}}}
\newcommand{\ba}{\begin{array}}
\newcommand{\ea}{\end{array}}
\newcommand{\on}{\operatorname}
\newcommand{\eq}{\begin{eqnarray}}
\newcommand{\eneq}{\end{eqnarray}}
\newcommand{\hs}{\hspace*}

\title[Decorated Geometric Crystals and Realizations of Crystal Bases]
{Decorated Geometric Crystals, 
Polyhedral and Monomial Realizations of Crystal Bases}

\author{Toshiki N\textsc{akashima}}
\address{Department of Mathematics, 
Sophia University, Kioicho 7-1, Chiyoda-ku, Tokyo 102-8554,
Japan}
\email{toshiki@sophia.ac.jp}
\thanks{supported in part by JSPS Grants 
in Aid for Scientific Research $\sharp 22540031$.}

\date{}

\dedicatory{Dedicated to Professor Michio Jimbo on the occasion of 
his 60th birthday}


\begin{abstract}
We shall show that for type $A_n$ 
the realization of crystal bases 
obtained from the decorated geometric crystals in \cite{BK2} coincides
 with our polyhedral realizations of crystal bases. We also observe
 certain relations of decorations and monomial realizations of crystal bases.
\end{abstract}

\maketitle
\renewcommand{\thesection}{\arabic{section}}
\section{Introduction}
\setcounter{equation}{0}
\renewcommand{\theequation}{\thesection.\arabic{equation}}

In \cite{BK2}, Berenstein and Kazhdan introduced the notion of decorated 
geometric crystals for reductive algebraic groups. 
Geometric crystals are geometric analogue to the
Kashiwara's crystal bases (\cite{BK}). We, indeed, treated geometric
crystals in the affine/Kac-Moody settings (\cite{KNO,KNO2,N}), but 
we do not need such general settings and then we shall consider the
(semi-)simple settings below.
Let $I$ be a finite index set. Associated with 
a Cartan matrix $A=(a_{i,j})_{i,j\in I}$, 
define the decorated geometric crystal $\cX=(\chi,f)$, which 
is a pair of geometric crystal 
$\chi=(X,\{e_i\}_i,\{\gamma_i\}_i,\{\vep_i\}_i)$ 
and  a certain special rational function $f$ such that 
\[
 f(e_i^c(x))=f(x)+(c-1)\vp_i(x)+(c^{-1}-1)\vep_i(x),
\]
for any $i\in I$, where $e_i^c$ is the rational $\bbC^\times$ action on
$X$, and $\vp_i$ and $\vep_i=\vep_i\cdot \gamma_i$ 
are the rational functions on $X$.

If we apply the procedure called ``ultra-discretization''(UD) to ``positive 
geometric crystals'' (see \ref{subsec-posi}), 
then we would obtain certain free-crystals for
the transposed Cartan matrix (\cite{BK,N}).  As for a positive decorated
geometric crystal $(\chi,f,T',\theta)$
applying  UD to the function $f$ and
considering the convex polyhedral domain defined by the inequality 
$UD(f)\geq0$, we get the crystal with the property
``normal''(\cite{K3}). Moreover, abstracting a connected component with 
the highest weight $\lm$, 
we obtain the Langlands dual 
Kashiwara's crystal $B(\lm)$ with the highest weight
$\lm$.

This result makes us recall the ``polyhedral realization'' of crystal
bases (\cite{N2,NZ}) since it has very similar way to get the crystal 
$B(\lm)$ from certain free-crystals, defined by the system of linear
inequalities. 
Thus, one of the main aims of this article is to show that the crystals
obtained by UD from positive decorated geometric crystals and the
polyhedral realizations of crystals coincide with each other for type $A_n$.

One more aim of this article is to describe the relations between 
the function $f_B$ for certain decorated geometric crystal
$(TB^-_{w_0},f_B)$
and monomial realization of crystals (\cite{K7,Nj}).
We shall propose the conjecture of their relations and 
present the affirmative answer for type $A_n$.
Let us mention the statement of the conjecture:
for the function $f_B$ and certain positive structure 
$T\Theta^-_{\bfi}$ on $TB^-_{w_0}$, the function 
$f_B(t\Theta^-_{\bf i}(c))$ is
expressed as a sum of monomials in the crystal $\cY(p)$ with positive
coefficients (for more details, see Conjecture \ref{conj1} below.).

Observing this relation, we can deduce the refined polyhedral
realization of crystals induced from the monomial realizations.
Indeed, for the original polyhedral realizations we are forced the
condition ``ample'', which is some technical condition to guarantee the
non-emptiness of the underlying crystal (see Theorem \ref{main}).
But, if the relations among the polyhedral realizations, the UD of decorated
geometric crystals and the monomial realizations are established, it would 
be possible to remove the condition ``ample'' and it would become easier
to obtain polyhedral realizations of crystals than applying
the present method.

The organization of the article is as follows:
in Sect.2, we review the theory of crystals and their polyhedral
realizations. In Sect.3, first we introduce the theory of decorated geometric 
crystals following \cite{BK2}.  Next, we define the decoration by using the 
elementary characters and certain special positive decorated geometric
crystal on ${\bbB}_w=TB^-_w$. Finally, the ultra-discretization of $TB^-_w$ is
described explicitly. We calculate the function $f_B$ exactly for type
$A_n$ in Sect.4. In Sect.5, for the type $A_n$ 
the coincidence of the polyhedral
realization $\Sigma_{\io}[\lm]$ and the ultra-discretization 
$B_{f_B,\Theta^-_{\bfii0}}(\lm)$ will be clarified by using the result in
Sect.4.
In the last section, we review the monomial realization of crystals
(\cite{K7,Nj})
and the function $f_B$ is expressed in terms of 
the monomials in the monomial realizations of crystals for type $A_n$. 
Finally, the conjecture is proposed and under the validity of the
conjecture,
we shall state the refined polyhedral realizations associated with the
monomial realizations.

The results for other simple Lie algebras are mentioned in the
forthcoming paper.

\renewcommand{\thesection}{\arabic{section}}
\section{Crystal and its polyhedral realization}
\setcounter{equation}{0}
\renewcommand{\theequation}{\thesection.\arabic{equation}}

\subsection{Notations}

We list the notations used in this paper.
Indeed, the settings below are originally Kac-Moody ones, but in the
article we do not need them and then we restrict the settings to
semi-simple ones.
Let $\ge$ be
a  semi-simple Lie algebra over $\bbQ$
with a Cartan subalgebra
$\tt$, a weight lattice $P \subset \tt^*$, the set of simple roots
$\{\al_i: i\in I\} \subset \tt^*$,
and the set of coroots $\{h_i: i\in I\} \subset \tt$,
where $I$ is a finite index set.
Let $\lan h,\lm\ran=\lm(h)$ be the pairing between $\tt$ and $\tt^*$,
and $(\al, \beta)$ be an inner product on
$\tt^*$ such that $(\al_i,\al_i)\in 2\bbZ_{\geq 0}$ and
$\lan h_i,\lm\ran={{2(\al_i,\lm)}\over{(\al_i,\al_i)}}$
for $\lm\in\tt^*$ and $A:=(\lan h_i,\al_j\ran)_{i,j}$ is the associated Cartan matrix.
Let $P^*=\{h\in \tt: \lan h,P\ran\subset\ZZ\}$ and
$P_+:=\{\lm\in P:\lan h_i,\lm\ran\in\ZZ_{\geq 0}\}$.
We call an element in $P_+$ a {\it dominant integral weight}.
The quantum algebra $\uq$
is an associative
$\QQ(q)$-algebra generated by the $e_i$, $f_i \,\, (i\in I)$,
and $q^h \,\, (h\in P^*)$
satisfying the usual relations.
The algebra $\uqm$ is the subalgebra of $\uq$ generated 
by the $f_i$ $(i\in I)$.

For the irreducible highest weight module of $\uq$
with the highest weight $\lm\in P_+$, we denote it by $V(\lm)$
and its {\it crystal base} we denote $(L(\lm),B(\lm))$.
Similarly, for the crystal base of the algebra $\uqm$ we denote 
$(L(\ify),B(\ify))$ (see \cite{K0,K1}).
Let $\pi_{\lm}:\uqm\longrightarrow V(\lm)\cong \uqm/{\sum_i\uqm
f_i^{1+\lan h_i,\lm\ran}}$ be the canonical projection and 
$\widehat \pi_{\lm}:L(\ify)/qL(\ify)\longrightarrow L(\lm)/qL(\lm)$
be the induced map from $\pi_{\lm}$. Here note that 
$\widehat \pi_{\lm}(B(\ify))=B(\lm)\sqcup\{0\}$.

By the terminology {\it crystal } we mean some combinatorial object 
obtained by abstracting the properties of crystal bases.
Indeed, crystal constitutes a set $B$ and the maps
$wt:B\longrightarrow P$, $\vep_i,\vp_i:B\longrightarrow \ZZ\sqcup\{-\ify\}$
and $\eit,\fit:B\sqcup\{0\}\longrightarrow B\sqcup\{0\}$
($i\in I$) satisfying several axioms (see \cite{K3},\cite{NZ},\cite{N2}).
In fact, $B(\ify)$ and $B(\lm)$ are the typical examples 
of crystals.

Let $B_1$ and $B_2$ be crystals.
A {\sl strict morphism } of crystals $\psi:B_1\lar B_2$
is a map $\psi:B_1\sqcup\{0\} \lar B_2\sqcup\{0\}$
satisfying the following conditions: $\psi(0)=0$, 
$wt(\psi(b)) = wt(b)$,
$\vep_i(\psi(b)) = \vep_i(b),$
$\vp_i(\psi(b)) = \vp_i(b)$
if $b\in B_1$ and $\psi(b)\in B_2,$
and the map $\psi: B_1\sqcup\{0\} \lar B_2\sqcup\{0\}$
commutes with all $\eit$ and $\fit$.
An injective strict morphism is called an {\it embedding }of crystals.

It is well-known that $\uq$ has a Hopf algebra structure.
Then the tensor product of $\uq$-modules has
a $\uq$-module structure.
The crystal bases have very nice properties for 
tensor operations. Indeed, if $(L_i,B_i)$ is a crystal base of 
$\uq$-module $M_i$ ($i=1,2$), $(L_1\ot_A L_2, B_1\ot B_2)$
is a crystal base of $M_1\ot_{\QQ(q)} M_2$ (\cite{K1}).
Consequently, we can consider the tensor product
of crystals and then they constitute a tensor category.

\subsection{Polyhedral Realization of $B(\ify)$}
\label{poly-uqm}
Let us recall the results in \cite{NZ}.

%
Consider the infinite $\bbZ$-lattice
\begin{equation}
\ZZ^{\ify}
:=\{(\cd,x_k,\cd,x_2,x_1): x_k\in\ZZ
\,\,{\rm and}\,\,x_k=0\,\,{\rm for}\,\,k\gg 0\};
\label{uni-cone}
\end{equation}
we will denote by $\ZZ^{\ify}_{\geq 0} \subset \ZZ^{\ify}$
the subsemigroup of nonnegative sequences.
To the rest of this section, we fix an infinite sequence of indices
$\io=\cd,i_k,\cd,i_2,i_1$ from $I$ such that
\begin{equation}
{\hbox{
$i_k\ne i_{k+1}$ and $\sharp\{k: i_k=i\}=\ify$ for any $i\in I$.}}
\label{seq-con}
\end{equation}

We can associate to $\io$ a crystal structure
on $\ZZ^{\ify}$ and denote it by $\ZZ^{\ify}_{\io}$ 
(\cite[2.4]{NZ}).

\begin{pro}[\cite{K3}, See also \cite{NZ}]
\label{emb}
There is a unique strict embedding of crystals
$($called Kashiwara embedding$)$
\begin{equation}
\Psi_{\io}:B(\ify)\hookrightarrow \ZZ^{\ify}_{\geq 0}
\subset \ZZ^{\ify}_{\io},
\label{psi}
\end{equation}
such that
$\Psi_{\io} (u_{\ify}) = (\cd,0,\cd,0,0)$, where 
$u_{\ify}\in B(\ify)$ is the vector corresponding to $1\in \uqm$.
\end{pro}

Consider the infinite dimensional vector space
$$
\QQ^{\ify}:=\{{x}=
(\cd,x_k,\cd,x_2,x_1): x_k \in \QQ\,\,{\rm and }\,\,
x_k = 0\,\,{\rm for}\,\, k \gg 0\},
$$
and its dual space $(\QQ^{\ify})^*:={\rm Hom}(\QQ^{\ify},\QQ)$.
We will write a linear form $\vp \in (\QQ^{\ify})^*$ as
$\vp({x})=\sum_{k \geq 1} \vp_k x_k$ ($\vp_j\in \QQ$)
for $x\in \QQ^{\ify}$.

For the fixed infinite sequence
$\io=(i_k)$ and $k\geq1$ we set $\kp:={\rm min}\{l:l>k\,\,{\rm and }\,\,i_k=i_l\}$ and
$\km:={\rm max}\{l:l<k\,\,{\rm and }\,\,i_k=i_l\}$ if it exists,
or $\km=0$  otherwise.
We set for $x\in \QQ^{\ify}$, $\beta_0(x)=0$ and
\begin{equation}
\beta_k(x):=x_k+\sum_{k<j<\kp}\lan h_{i_k},\al_{i_j}\ran x_j+x_{\kp}
\qq(k\geq1).
\label{betak}
\end{equation}
We define the piecewise-linear operator 
$S_k=S_{k,\io}$ on $(\QQ^{\ify})^*$ by
$$
S_k(\vp):=
\left\{
\begin{array}{ll}
\vp-\vp_k\beta_k & {\mbox{ if }}\vp_k>0,\\
 \vp-\vp_k\beta_{\km} & {\mbox{ if }}\vp_k\leq 0.
\end{array}
\right.
$$
Here we set
\begin{eqnarray}
\Xi_{\io} &:=  &\{S_{j_l}\cd S_{j_2}S_{j_1}x_{j_0}\,|\,
l\geq0,j_0,j_1,\cd,j_l\geq1\},
\label{Xi_io}\\
\Sigma_{\io} & := &
\{x\in \ZZ^{\ify}\subset \QQ^{\ify}\,|\,\vp(x)\geq0\,\,{\rm for}\,\,
{\rm any}\,\,\vp\in \Xi_{\io}\}.
\end{eqnarray}
We impose on $\io$ the following positivity assumption:
\begin{equation}
{\hbox{if $\km=0$ then $\vp_k\geq0$ for any 
$\vp(x)=\sum_k\vp_kx_k\in \Xi_{\io}$}}.
\label{posi}
\end{equation}
\begin{thm}[\cite{NZ}]
Let $\io$ be a sequence of indices satisfying $(\ref{seq-con})$ 
and (\ref{posi}). Then we have 
${\rm Im}(\Psi_{\io})(\cong B(\ify))=\Sigma_{\io}$.
\end{thm}

\subsection{Structure of $\ZZ^{\ify}_{\io}[\lm]$}
\label{ZZ-io-lm}

Let $R_{\lm}:=\{r_{\lm}\}$ be the crystal 
which consists of one element $r_{\lm}$ (\cite{N2}).
Consider the crystal $\ZZ^{\ify}_{\io}\ot R_{\lm}$ 
and denote it by 
$\ZZ^{\ify}_{\io}[\lm]$. Here note that
since the crystal $R_{\lm}$ has  only one element,
as a set we can identify $\ZZ^{\ify}_{\io}[\lm]$ with
$\ZZ^{\ify}_{\io}$ but their crystal structures are different.
So we review an explicit crystal structure of
$\ZZ^{\ify}[\lm]$ in \cite{N2}.
Fix a sequence of indices $\io:=(i_k)_{k\geq 1}$ satisfying the condition
(\ref{seq-con}) and a weight $\lm\in P$
(Here we do not necessarily assume that 
$\lm$ is dominant.).
$\ZZ^{\ify}[\lm]$ can be regarded 
as a subset of $\QQ^{\ify}$, and then 
we denote an element in $\ZZ^{\ify}[\lm]$
by $ x=(\cd,x_k,\cd,x_2,x_1)$.
For $ x=(\cd,x_k,\cd,x_2,x_1)\in \QQ^{\ify}$
we define the linear functions
\begin{eqnarray}
\sigma_k(x)&:= &x_k+\sum_{j>k}\lan h_{i_k},\al_{i_j}\ran x_j,
\q(k\geq1)
\label{sigma}\\
\sigma_0^{(i)}(x)
&:= &-\lan h_i,\lm\ran+\sum_{j\geq1}\lan h_i,\al_{i_j}\ran x_j,
\q(i\in I)
\label{sigma0}
\end{eqnarray}
Here note that
since $x_j=0$ for $j\gg0$ on $\QQ^{\ify}$,
the functions $\sigma_k$ and $\sigma^{(i)}_0$ are
well-defined.
Let $\sigma^{(i)} (x)
 := {\rm max}_{k: i_k = i}\sigma_k (x)$, and
$M^{(i)} :=
\{k: i_k = i, \sigma_k (x) = \sigma^{(i)}(x)\}.
$
Note that
$\sigma^{(i)} (x)\geq 0$, and that
$M^{(i)} = M^{(i)} (x)$ is a finite set
if and only if $\sigma^{(i)} (x) > 0$.
Now we define the maps
$\eit: \ZZ^{\ify}[\lm] \sqcup\{0\}\lar \ZZ^{\ify}[\lm] \sqcup\{0\}$
and
$\fit: \ZZ^{\ify}[\lm] \sqcup\{0\}\lar \ZZ^{\ify}[\lm] \sqcup\{0\}$ 
by setting $\eit(0)=\fit(0)=0$, and 
\begin{equation}
(\fit(x))_k  = x_k + \delta_{k,{\rm min}\,M^{(i)}}
\,\,{\rm if }\,\,\sigma^{(i)}(x)>\sigma^{(i)}_0(x);
\,\,{\rm otherwise}\,\,\fit(x)=0,
\label{action-f}
\end{equation}
\begin{equation}
(\eit(x))_k  = x_k - \delta_{k,{\rm max}\,M^{(i)}} \,\, {\rm if}\,\,
\sigma^{(i)} (x) > 0\,\,
{\rm and}\,\,\sigma^{(i)}(x)\geq\sigma^{(i)}_0(x) ; \,\,
 {\rm otherwise} \,\, \eit(x)=0,
\label{action-e}
\end{equation}
where $\del_{i,j}$ is the Kronecker's delta.
We also define the functions
$wt$, $\vep_i$ and $\vp_i$ on $\ZZ^{\ify}[\lm]$ by
\begin{eqnarray}
&& wt(x) :=\lm -\sum_{j=1}^{\ify} x_j \al_{i_j},
\label{wt-vep-vp-1}\\
&& \vep_i (x) := {\rm max}(\sigma^{(i)} (x),\sigma_0^{(i)}(x))
\label{wt-vep-vp}\\
&& \vp_i (x) := \lan h_i, wt(x) \ran + \vep_i(x).
\label{wt-vep-vp-3}
\end{eqnarray}

Note that 
by (\ref{wt-vep-vp-1}) we have
\begin{equation}
\lan h_i,wt(x)\ran = -\sigma^{(i)}_0(x).
\label{**}
\end{equation}

\subsection{Polyhedral Realization of $B(\lm)$}

In this subsection, 
we review the result in \cite{N2}. 
In the rest of this section,
$\lm$ is supposed to be a dominant integral weight.
Here we define the map
\begin{equation}
\Phi_{\lm}:(B(\ify)\ot R_{\lm})\sqcup\{0\}\longrightarrow B(\lm)\sqcup\{0\},
\label{philm}
\end{equation}
by $\Phi_{\lm}(0)=0$ and $\Phi_{\lm}(b\ot r_{\lm})=\wpi(b)$ for $b\in B(\ify)$.
We set
$$
\wtil B(\lm):=
\{b\ot r_{\lm}\in B(\ify)\ot R_{\lm}\,|\,\Phi_{\lm}(b\ot r_{\lm})\ne 0\}.
$$

\begin{thm}[\cite{N2}]
\label{ify-lm}
\begin{enumerate}
\item
The map $\Phi_{\lm}$ becomes a surjective strict morphism of crystals
$B(\ify)\ot R_{\lm}\longrightarrow B(\lm)$.
\item
$\wtil B(\lm)$ is a subcrystal of $B(\ify)\ot R_{\lm}$, 
and $\Phi_{\lm}$ induces the
isomorphism of crystals $\wtil B(\lm)\mapright{\sim} B(\lm)$.
\end{enumerate}
\end{thm}

By Theorem \ref{ify-lm}, we have the strict embedding of crystals
$\Omega_{\lm}:B(\lm)(\cong \wtil B(\lm))\hookrightarrow B(\ify)\ot R_{\lm}.$
Combining $\Omega_{\lm}$ and the
Kashiwara embedding $\Psi_{\io}$,
we obtain the following:
\begin{thm}[\cite{N2}]
\label{embedding}
There exists the unique  strict embedding of crystals
\begin{equation}
\Psi_{\io}^{(\lm)}:B(\lm)\stackrel{\Omega_{\lm}}{\hookrightarrow}
B(\ify)\ot R_{\lm}
\stackrel{\Psi_{\io}\ot {\rm id}}{\hookrightarrow}
\ZZ^{\ify}_{\io}[\lm],
\label{Psi-lm}
\end{equation}
such that $\Psi^{(\lm)}_{\io}(u_{\lm})=(\cd,0,0,0)\ot r_{\lm}$.
\end{thm}

\vskip5pt


We fix a sequence of indices $\io$
satisfying (\ref{seq-con}) and take a dominant integral weight 
$\lm\in P_+$.
For $k\geq1$ let $k^{(\pm)}$ be  the ones in \ref{poly-uqm}.
Let $\beta_k^{(\pm)}(x)$ be linear functions given by
\begin{eqnarray}
&& \q
\beta_k^{(+)} (x)  =  \sigma_k (x) - \sigma_{\kp} (x)
= x_k+\sum_{k<j<\kp}\lan h_{i_k},\al_{i_j}\ran x_j+x_{\kp},
\label{beta}\\
&&  \beta_k^{(-)} (x) 
=
\left\{
\begin{array}{ll}
\sigma_{\km} (x) - \sigma_k (x)
=x_{\km}+\sum_{\km<j<k}\lan h_{i_k},\al_{i_j}\ran x_j+x_k &
 \hspace{-10pt} {\mbox{ if }}\km>0,\\
\sigma_0^{(i_k)} (x) - \sigma_k (x)
=-\lan h_{i_k},\lm\ran+\sum_{1\leq j<k}\lan h_{i_k},\al_{i_j}\ran x_j+x_k
&
 \hspace{-10pt}{\mbox{ if }}\km=0,
\end{array}
\right.\label{beta--}
\end{eqnarray}

Here note that
$\beta_k^{(+)}=\beta_k$ and 
$\beta_k^{(-)}=\beta_{\km}  {\hbox{ \,\,if\,\, $\km>0$}}$.

Using this notation, for every $k \geq 1$, we define 
an operator
$\what S_k = \what S_{k,\io}$ for a linear function 
$\vp(x)=c+\sum_{k\geq 1}\vp_kx_k$
$(c,\vp_k\in\QQ)$ on $\QQ^{\ify}$ by:

$$
\what S_k\,(\vp) :=\left\{
\begin{array}{ll}
\vp - \vp_k \beta_k^{(+)} & {\mbox{ if }}\vp_k > 0,\\
\vp - \vp_k \beta_k^{(-)} & {\mbox{ if }}\vp_k \leq 0.
\end{array}
\right.
$$

For the fixed sequence $\io=(i_k)$, 
in case $\km=0$ for $k\geq1$, there exists unique $i\in I$ such that $i_k=i$.
We denote such $k$ by $\io^{(i)}$, namely, $\io^{(i)}$ is the first
number $k$ such that $i_k=i$.
Here for $\lm\in P_+$ and $i\in I$ we set
\begin{equation}
\lm^{(i)}(x):=
-\beta^{(-)}_{\io^{(i)}}(x)=\lan h_i,\lm\ran-\sum_{1\leq j<\io^{(i)}}
\lan h_i,\al_{i_j}\ran x_j-x_{\io^{(i)}}.
\label{lmi}
\end{equation}

For $\io$ and a dominant integral weight $\lm$,
let $\Xi_{\io}[\lm]$ be the set of all linear functions
generated by $\what S_k=\what S_{k,\io}$ 
from the functions $x_j$ ($j\geq1$)
and $\lm^{(i)}$ ($i\in I$), namely,
\begin{equation}
\begin{array}{ll}
\Xi_{\io}[\lm]&:=\{\what S_{j_l}\cd\what S_{j_1}x_{j_0}\,
:\,l\geq0,\,j_0,\cd,j_l\geq1\}
\\
&\cup\{\what S_{j_k}\cd \what S_{j_1}\lm^{(i)}(x)\,
:\,k\geq0,\,i\in I,\,j_1,\cd,j_k\geq1\}.
\end{array}
\label{Xi}
\end{equation}
Now we set
\begin{equation}
\Sigma_{\io}[\lm]
:=\{x\in \ZZ^{\ify}_{\io}[\lm](\subset \QQ^{\ify})\,:\,
\vp(x)\geq 0\,\,{\rm for \,\,any }\,\,\vp\in \Xi_{\io}[\lm]\}.
\label{Sigma}
\end{equation}

For a sequence $\io$ and a dominant integral weight $\lm$, a pair
$(\io,\lm)$ is called {\it ample}
if $\Sigma_{\io}[\lm]\ni\vec 0=(\cd,0,0)$.

\begin{thm}[\cite{N2}]
\label{main}
Suppose that $(\io,\lm)$ is ample.
Then we have
${\rm Im}(\plm)(\cong B(\lm))=\Sigma_{\io}[\lm]$,
where the explicit form of $\vep_i$ on $\Sigma_{\io}[\lm]$ is as follows:
\begin{equation}
\vep_i(x) =\sigma^{(i)}(x).
\label{wt-vep-vp4}
\end{equation}
The other formula for $\vp_i$, $\eit$ and $\fit$ are same as above.
\end{thm}
{\sl Proof.}
The formula (\ref{wt-vep-vp4})  slightly differs from (\ref{wt-vep-vp}).
Indeed, by (\ref{action-e}) we know that for $x\in \Sigma_{\io}[\lm]$
unless $\sigma^{(i)}(x)>0$ and $\sigma^{(i)}(x)\geq \sigma^{(i)}_0(x)$, 
we find $\eit(x)=0$. Furthermore, for any $x=(\cd,x_2,x_1)\in
\Sigma_\io[\lm]$
it follows from the definition of
$\Sigma_\io[\lm]$ that $0\leq \lm^{(i)}(x)
=\sigma_{\io^{(i)}}(x)-\sigma_0^{(i)}(x)$
which implies that $\sigma^{(i)}(x)\geq \sigma_0^{(i)}(x)$ and then
we can obtain (\ref{wt-vep-vp4}).\qed

\medskip
\subsection{\bf $A_n$-case}
We shall apply the results in the previous subsection to 
the case $\ge=A_n$. 
Let us identify the index set $I$ with $[1,n] := \{1,2,\cd,n\}$ 
in the standard way; thus, the Cartan matrix 
$(a_{i,j}= \lan h_i,\al_j\ran )_{1 \leq i,j \leq n}$ is given by 
$a_{i,i}=2$, $a_{i,j}=-1$ for $|i-j|=1$, and 
$a_{i,j}=0$ otherwise. 
As the infinite sequence $\io$ let us take 
the following periodic sequence 
$$
\io = \cd,\underbrace{n,\cd,2,1}_{},
\cd,\underbrace{n,\cd,2,1}_{},\underbrace{n,\cd,2,1}_{}.
$$
Following \cite[Sect.5]{NZ}, we shall 
change the indexing set for $\ZZ^{\ify}$
from $\ZZ_{\geq 1}$ to $\ZZ_{\geq 1} \times [1,n]$, which is given by 
the bijection 
$\ZZ_{\geq 1} \times [1,n] \to \ZZ_{\geq 1}$ 
($(j;i) \mapsto (j-1)n + i$). 
According to this, we will write an element $x \in \ZZ^{\ify}$
as a doubly-indexed family $(x_{j;i})_{j \geq 1, i \in [1,n]}$.
We will adopt the convention that $x_{j;i} = 0$ unless
$j \geq 1$ and  $i \in [1,n]$; in particular, $x_{j;0} = x_{j;n+1} = 0$
for all $j$. 

\begin{thm} 
\label{A_n}
Let $\lm=\sum_{1\leq i\leq n}\lm_i\Lm_i$ $(\lm_i\in \ZZ_{\geq0})$ 
be a dominant integral weight. 
In the above notation, the image ${\rm Im} \,(\Psi^{(\lm)}_{\io})$ 
is the set of all integer families $(x_{j;i})$ such that 
\beqn
&&\hspace{-30pt}\hbox{
$x_{1;i} \geq x_{2;i-1} \geq \cd \geq x_{i;1} \geq 0$ 
for $1 \leq i \leq n$}
\label{sl-1}\\
&&\hspace{-30pt}\hbox{
$x_{j;i} = 0$ for $i+j > n+1$, }
\label{j;i=0}\\
&&\hspace{-30pt}\hbox{
$\lm_i\geq x_{j;i-j+1}-x_{j;i-j}$ for
$1\leq j\leq i\leq n$.}
\label{sl-2}
\eeqn
\end{thm}
Observing (\ref{j;i=0}), we can rewrite the theorem in the following form:
Let $\io_0$ be one of the reduced longest words of type $A_n$:
\begin{equation}
\io_0=\underbrace{1}_{},
\underbrace{2,1}_{},\underbrace{3,2,1}_{},\cd 
\underbrace{n,n-1,\cd,2,1}_{}.
\label{io0}
\end{equation}
\begin{cor}
\label{cor-A_n}
Associated with $\io_0$, we define 
\begin{equation}
\bbZ_{\io_0}[\lm]:=\{(x_{j;i}|1\leq i+j\leq n+1)\in
 \bbZ^{\frac{(n(n+1))}{2}}|
(x_{j;i})\text{ satisfies (\ref{sl-1}) and (\ref{sl-2}).}\}
\end{equation}
There exists the crystal structure on $\bbZ_{\io_0}[\lm]$ induced
from the one on $\bbZ_\io[\lm]$ and then the crystal $\bbZ_{\io_0}[\lm]$
is isomorphic to $B(\lm)$.
\end{cor}

\renewcommand{\thesection}{\arabic{section}}
\section{Decorated geometric crystals}
\setcounter{equation}{0}
\renewcommand{\theequation}{\thesection.\arabic{equation}}

The basic reference for this section is \cite{BK,BK2}.
\subsection{Definitions}

Let  
 $A=(a_{ij})_{i,j\in I}$ be 
an indecomposable Cartan matrix
with a finite index set $I$
(though we can consider more general Kac-Moody setting.).
Let $(\frt,\{\al_i\}_{i\in I},\{h_i\}_{i\in I})$ 
be the associated
root data 
satisfying $\al_j(h_i)=a_{ij}$.
Let $\ge=\ge(A)=\lan \frt,e_i,f_i(i\in I)\ran$ 
be the simple Lie algebra associated with $A$
over $\bbC$ and $\Del=\Del_+\sqcup\Del_-$
be the root system associated with $\ge$, where $\Del_{\pm}$ is 
the set of positive/negative roots.

Define the simple reflections $s_i\in{\rm Aut}(\frt)$ $(i\in I)$ by
$s_i(h)\seteq h-\al_i(h)h_i$, which generate the Weyl group $W$.
Let $G$ be the simply connected simple algebraic group 
over $\bbC$ whose Lie algebra is $\ge=\frn_+\oplus \frt\oplus \frn_-$, 
which is the usual triangular decomposition.
Let $U_{\al}\seteq\exp\ge_{\al}$ $(\al\in \Del)$
be the one-parameter subgroup of $G$.
The group $G$ (resp. $U^\pm$) is generated by 
$\{U_{\al}|\al\in \Del\}$ 
(resp. $\{U_\al|\al\in\Del_{\pm}$).
Here $U^\pm$ is a unipotent radical of $G$ and 
${\rm Lie}(U^\pm)=\frn_{\pm}$.
For any $i\in I$, there exists 
a unique group homomorphism
$\phi_i\cl SL_2(\bbC)\rightarrow G$ such that
\[
\phi_i\left(
\left(
\begin{array}{cc}
1&t\\
0&1
\end{array}
\right)\right)=\exp(t e_i),\,\,
 \phi_i\left(
\left(
\begin{array}{cc}
1&0\\
t&1
\end{array}
\right)\right)=\exp(t f_i)\qquad(t\in\bbC).
\]
Set $\al^\vee_i(c)\seteq
\phi_i\left(\left(
\begin{smallmatrix}
c&0\\
0&c^{-1}\end{smallmatrix}\right)\right)$,
$x_i(t)\seteq\exp{(t e_i)}$, $y_i(t)\seteq\exp{(t f_i)}$, 
$G_i\seteq\phi_i(SL_2(\bbC))$,
$T_i\seteq \alpha_i^\vee(\bbC^\times)$ 
and 
$N_i\seteq N_{G_i}(T_i)$. Let
$T$ be a maximal torus of $G$ 
which has $P$ as its weight lattice and 
Lie$(T)=\frt$.
Let 
$B^{\pm}(\supset T)$ be the Borel subgroup of $G$.
We have the isomorphism
$\phi:W\mapright{\sim}N/T$ defined by $\phi(s_i)=N_iT/T$.
An element $\ovl s_i:=x_i(-1)y_i(1)x_i(-1)$ is in 
$N_G(T)$, which is a representative of 
$s_i\in W=N_G(T)/T$.

\begin{df}
\label{def-gc}
Let $X$ be an affine algebraic variety over $\bbC$, 
$\gamma_i$, $\vep_i, f$ 
$(i\in I)$ rational functions on $X$, and 
$e_i:\bbC^\times\times X\to X$ a unital rational $\bbC^\times$-action.
A 5-tuple $\chi=(X,\{e_i\}_{i\in I},\{\gamma_i,\}_{i\in I},
\{\vep_i\}_{i\in I},f)$ is a 
$G$ (or $\ge$)-{\it decorated geometric crystal} 
if
\begin{enumerate}
\item
$(\{1\}\times X)\cap dom(e_i)$ 
is open dense in $\{1\}\times X$ for any $i\in I$, where
$dom(e_i)$ is the domain of definition of
$e_i\cl\C^\times\times X\to X$.
\item
The rational functions  $\{\gamma_i\}_{i\in I}$ satisfy
$\gamma_j(e^c_i(x))=c^{a_{ij}}\gamma_j(x)$ for any $i,j\in I$.
\item
The function $f$ satisfies
\begin{equation}
f(e_i^c(x))=f(x)+{(c-1)\vp_i(x)}+{(c^{-1}-1)\vep_i(x)},
\label{f}
\end{equation}
for any $i\in I$ and $x\in X$,  where $\vp_i:=\vep_i\cdot\gamma_i$.
\item
$e_i$ and $e_j$ satisfy the following relations:
\[
 \begin{array}{lll}
&\hspace{-20pt}e^{c_1}_{i}e^{c_2}_{j}
=e^{c_2}_{j}e^{c_1}_{i}&
{\rm if }\,\,a_{ij}=a_{ji}=0,\\
&\hspace{-20pt} e^{c_1}_{i}e^{c_1c_2}_{j}e^{c_2}_{i}
=e^{c_2}_{j}e^{c_1c_2}_{i}e^{c_1}_{j}&
{\rm if }\,\,a_{ij}=a_{ji}=-1,\\
&\hspace{-20pt}
e^{c_1}_{i}e^{c^2_1c_2}_{j}e^{c_1c_2}_{i}e^{c_2}_{j}
=e^{c_2}_{j}e^{c_1c_2}_{i}e^{c^2_1c_2}_{j}e^{c_1}_{i}&
{\rm if }\,\,a_{ij}=-2,\,
a_{ji}=-1,\\
&\hspace{-20pt}
e^{c_1}_{i}e^{c^3_1c_2}_{j}e^{c^2_1c_2}_{i}
e^{c^3_1c^2_2}_{j}e^{c_1c_2}_{i}e^{c_2}_{j}
=e^{c_2}_{j}e^{c_1c_2}_{i}e^{c^3_1c^2_2}_{j}e^{c^2_1c_2}_{i}
e^{c^3_1c_2}_je^{c_1}_i&
{\rm if }\,\,a_{ij}=-3,\,
a_{ji}=-1.
\end{array}
\]
\item
The rational functions $\{\vep_i\}_{i\in I}$ satisfy
$\vep_i(e_i^c(x))=c^{-1}\vep_i(x)$ and 
$\vep_i(e_j^c(x))=\vep_i(x)$ if $a_{i,j}=a_{j,i}=0$.
\end{enumerate}
\end{df}
We call the function $f$ in (iii) the {\it decoration} of $\chi$ and 
the relations in (iv) are called 
{\it Verma relations}.
If $\chi=(X,\{e_i\},\,\{\gamma_i\},\{\vep_i\})$
satisfies the conditions (i), (ii), (iv) and (v), 
we call $\chi$ a {\it geometric crystal}.
{\sl Remark.}
The definitions of $\vep_i$ and $\vp_i$ are different from the ones in 
e.g., \cite{BK2} since we adopt the definitions following
\cite{KNO,KNO2}. Indeed, if we flip $\vep_i\to \vep^{-1}$ and 
$\vp_i\to \vp^{-1}$, they coincide with ours.
\subsection{Characters}

Let $\what U:={\rm Hom}(U,\bbC)$ be the set of additive characters 
of $U$.
The {\it elementary character }$\chi_i\in \what U$ and  
the {\it standard regular character} $\chi^{\rm st}\in \what U$ are  
defined by 
\[
\chi_i(x_j(c))=\del_{i,j}\cdot c \q(c\in \bbC,\,\, i\in I),\qq
\chi^{st}=\sum_{i\in I}\chi_i.
\]
Let us define an anti-automorphism $\eta:G\to G$ 
by 
\[
 \eta(x_i(c))=x_i(c),\q  \eta(y_i(c))=y_i(c),\q \eta(t)=t^{-1}\q
 (c\in\bbC,\,\, t\in T),
\]
which is called the {\it positive inverse}.

The rational function $f_B$ on $G$ is defined by 
\begin{equation}
f_B(g)=\chi^{st}(\pi^+(w_0^{-1}g))+\chi^{st}(\pi^+(w_0^{-1}\eta(g))),
\label{f_B}
\end{equation}
for $g\in B\ovl w_0 B$, where $\pi^+:B^-U\to U$ is the projection by 
$\pi^+(bu)=u$.

For a split algebraic torus $T$ over $\bbC$, let us denote 
its lattice of (multiplicative )characters(resp. co-characters) by $X^*(T)$ 
(resp. $X_*(T)$). By the usual way, we identify $X^*(T)$
(resp. $X_*(T)$) with the weight lattice $P$ (resp. the dual weight
lattice $P^*$). 

\subsection{Positive structure and ultra-discretization}
\label{subsec-posi}

In this subsection, we review the notion positive structure and 
the ultra-discretization, which is called the 
tropicalization in \cite{BK,BK2}.

\begin{df}
Let $T,T'$ be  split algebraic tori over $\bbC$.
\begin{enumerate}
\item
A regular function $f=\sum_{\mu\in X^*(T)}c_\mu\cdot\mu$ on $T$ 
is {\it positive} if all coefficients $c_\mu$ are non-negative 
numbers. A rational function on $T$ is said to be {\it positive} if 
there exist positive regular functions $g,h$ such that 
$f=\frac{g}{h}$ ($h\ne0$).
\item
Let $f:T\to T'$ be a rational map between 
$T$ and $T'$. Then we say that $f$ is {\it positive} if 
for any $\xi\in X^*(T')$ we have that $\xi\circ f$ is positive in the
above sense.
\end{enumerate}
\end{df}
Note that if $f,g$ are positive rational functions on $T$, then 
$f\cdot g$, $f/g$ and $f+g$ are all positive.

\begin{df}
Let $\chi=(X,\{e_i\}_{i\in I},\{{\rm wt}_i\}_{i\in I},
\{\vep_i\}_{i\in I},f)$ be a decorated
geometric crystal, $T'$ an algebraic torus
and $\theta:T'\rightarrow X$ 
a birational map.
The birational map $\theta$ is called 
{\it positive structure} on
$\chi$ if it satisfies:
\begin{enumerate}
\item For any $i\in I$ the rational functions
$\gamma_i\circ \theta, \vep_i\circ \theta, f\circ\theta:T'\rightarrow \bbC$ 
are all positive in the above sense.
\item
For any $i\in I$, the rational map
$e_{i,\theta}:\bbC^\tm \tm T'\rightarrow T'$ defined by
$e_{i,\theta}(c,t)
:=\theta^{-1}\circ e_i^c\circ \theta(t)$
is positive.
\end{enumerate}
\end{df}

Let $v:\bbC(c)\setminus\to\bbZ$ be a map defined by 
$v(f(c)):=\deg(f(c^{-1}))$,which is different from that in e.g.,
\cite{KNO,KNO2,N,N3}. Note that this definition of the map $UD$ is
called tropicalization in \cite{BK} and much simpler than the one 
in \cite{BK2} since it is sufficient in this article.
Here, we have the formula for positive rational functions $f$ and $g$:
\begin{equation}
v(f\cdot g)=v(f)+v(g),\q
v(f/g)=v(f)-v(g),\q
v(f+g)=\min(v(f),v(g)).
\label{uf-formula}
\end{equation}

Let $f\cl T\rightarrow T'$ be a 
positive rational mapping
of algebraic tori $T$ and 
$T'$.
We define a map $\what f\cl X_*(T)\rightarrow X_*(T')$ by 
\[
\langle\chi,\what f(\xi)\rangle
=v(\chi\circ f\circ \xi),
\]
where $\chi\in X^*(T')$ and $\xi\in X_*(T)$.

Let $\cT_+$ be the category whose objects are algebraic tori over
$\bbC$ and whose morphisms are positive rational maps.
Then, we obtain the functor
\[
\begin{array}{cccc}
{\mathcal UD}:& \cT_+&\longrightarrow &{\mathfrak Set}\\
&T& \mapsto &X_*(T)\\
&(f:T\to T') &\mapsto &(\what f:X_*(T)\to X_*(T')).
\end{array}
\]

Let $\theta:T\rightarrow X$ be a positive structure on 
a decorated geometric crystal $\chi=(X,\{e_i\}_{i\in I},
\{{\rm wt}_i\}_{i\in I},
\{\vep_i\}_{i\in I,f})$.
Applying the functor ${\mathcal UD}$ 
to positive rational morphisms
$e_{i,\theta}:\bbC^\tm \tm T'\rightarrow T'$ and
$\gamma\circ \theta:T'\ra T$
(the notations are
as above), we obtain
\begin{eqnarray*}
\til e_i&:=&{\mathcal UD}(e_{i,\theta}):
\ZZ\tm X_*(T) \rightarrow X_*(T)\\
{\rm wt}_i&:=&{\mathcal UD}(\gamma_i\circ\theta):
X_*(T')\rightarrow \bbZ,\\
\wtil\vep_i&:=&{\mathcal UD}(\vep_i\circ\theta):
X_*(T')\rightarrow \bbZ,\\
\wtil f&:=& {\mathcal UD}(f \circ\theta):
\end{eqnarray*}
Now, for given positive structure $\theta:T'\rightarrow X$
on a geometric crystal 
$\chi=(X,\{e_i\}_{i\in I},\{{\rm wt}_i\}_{i\in I},
\{\vep_i\}_{i\in I})$, we associate 
the quadruple $(X_*(T'),\{\til e_i\}_{i\in I},
\{{\rm wt}_i\}_{i\in I},\{\wtil\vep_i\}_{i\in I})$
with a free pre-crystal structure (see \cite[2.2]{BK}) 
and denote it by ${\mathcal UD}_{\theta,T'}(\chi)$.
We have the following theorem:

\begin{thm}[\cite{BK,BK2,N}]
For any geometric crystal 
$\chi=(X,\{e_i\}_{i\in I},\{\gamma_i\}_{i\in I},
\{\vep_i\}_{i\in I})$ and positive structure
$\theta:T'\rightarrow X$, the associated pre-crystal 
${\mathcal UD}_{\theta,T'}(\chi)=
(X_*(T'),\{e_i\}_{i\in I},\{{\rm wt}_i\}_{i\in I},
\{\wtil\vep_i\}_{i\in I})$ 
is a Kashiwara's crystal.
\end{thm}
{\sl Remark.}
The definition of $\wtil\vep_i$ is different from the one in 
\cite[6.1.]{BK2} since our definition of $\vep_i$ corresponds to 
$\vep_i^{-1}$ in \cite{BK2}.

Now, for a positive decorated geometric crystal 
$\cX=((X,\{e_i\}_{i\in I},\{\gamma_i\}_{i\in I},
\{\vep_i\}_{i\in I},f),\theta,T')$, set 
\begin{equation}
 \wtil B_{\wtil f}:=\{\wtil x\in X_*(T')|\wtil f(\wtil x)\geq0\},
\label{btil}
\end{equation}
where $X_*(T')$ is identified with $\bbZ^{\dim(T')}$. Define 
\begin{equation}
B_{f,\theta}:=(\wtil B_{\wtil f},\wt_i|_{\wtil B_{\wtil f}},
\vep_i|_{\wtil B_{\wtil f}},e_i|_{\wtil B_{\wtil f}})_{i\in I}.
\label{btheta}
\end{equation}
\begin{pro}[\cite{BK2}]
For a positive decorated geometric crystal 
$\cX=((X,\{e_i\}_{i\in I},\{\gamma_i\}_{i\in I},
\{\vep_i\}_{i\in I},f),\theta,T')$, the
 quadruple $B_{f,\theta}$ in (\ref{btheta}) is a normal crystal.
\end{pro}

\subsection{Decorated geometric crystal on $\bbB_w$}
For a Weyl group element $w\in W$, define $B^-_w$ by 
\begin{equation}
B^-_w:=B^-\cap U\ovl w U.
\label{B-w}
\end{equation}
Now, set  $\bbB_w:=TB^-_w$. 
Let $\gamma_i:\bbB_w\to\bbC$ be the rational function defined by 
\begin{equation}
\gamma_i:\bbB_w\,\,\hookrightarrow \,\,\,
B^-\,\,\mapright{\sim}\,\,T\times U^-\,\,
\mapright{\rm proj}\,\,\, T\,\,\,\mapright{\al_i^\vee}\,\,\,\bbC.
\label{gammai}
\end{equation}

For any $i\in I$, there exists the natural projection 
$pr_i:B^-\to B^-\cap \phi(SL_2)$. Hence, 
for any $x\in \bbB_w$ there exists unique
       $v=\begin{pmatrix}b_{11}&0\\b_{21}&b_{22}\end{pmatrix}
\in SL_2$ such that 
$pr_i(x)=\phi_i(v)$. Using this fact, we define 
the rational function $\vep_i$ on $\bbB_w$ by 
\begin{equation}
\vep_i(x)=\frac{b_{22}}{b_{21}}\q(x\in\bbB_w).
\label{vepi}
\end{equation}
The rational $\bbC^\times$-action $e_i$ on $\bbB_w$ is defined by
\begin{equation}
e_i^c(x):=x_i\left((c-1)\vp_i(x)\right)\cdot x\cdot
x_i\left((c^{-1}-1)\vep_i(x)\right)\qq
(c\in\bbC^\times,\,\,x\in \bbB_w),
\label{ei-action}
\end{equation}
if $\vep_i(x)$ is well-defined, that is, $b_{21}\ne0$, 
and $e_i^c(x)=x$ if $b_{21}=0$.\\
{\sl Remark.} The definition (\ref{vepi}) is different from the one in 
\cite{BK2}. Indeed, if we take $\vep_i(x)=b_{21}/b_{22}$, 
then it coincides with
the one in \cite{BK2}.
\begin{pro}[\cite{BK2}]
For any $w\in W$,
the 5-tuple $\chi:=(\bbB_w,\{e_i\}_i,\{\gamma_i\}_i,\{\vep_i\}_i,f_B)$
is a decorated geometric crystal, where 
$f_B$ is in (\ref{f_B}), $\gamma_i$ is in (\ref{gammai}), $\vep_i$ is in 
(\ref{vepi}) and $e_i$ is in (\ref{ei-action}).
\end{pro}

\def\ld{\ldots}
For the  longest Weyl group element $w_0\in W$, let 
$\bfii0=i_1\ld i_N$ be one of its reduced expressions and 
define the positive structure on $B^-_{w_0}$ 
$\Theta^-_\bfii0:(\bbC^\times)^N\longrightarrow B^-_{w_0}$ by 
\[
 \Theta^-_\bfii0(c_1,\cd,c_N):=\pmby_{i_1}(c_1)\cd \fry_{i_N}(c_N),
\]
where $\pmby_i(c)=y_i(c)\al^\vee(c^{-1})$, 
which is different from $Y_i(c)$ in 
\cite{N,N2,KNO,KNO2}. Indeed, $Y_i(c)=\pmby_i(c^{-1})$.
We also define the positive structure on $\bbB_{w_0}$ as
$T\Theta^-_\bfii0:T\times(\bbC^\times)^N\,\,\longrightarrow\,\,\bbB_{w_0}$  
by $T\Theta^-_\bfii0(t,c_1,\cd,c_N)
=t\Theta^-_\bfii0(c_1,\cd,c_N)$.

Now, for this positive structure, we describe the geometric crystal
structure on $\bbB_{w_0}=TB^-_{w_0}$ explicitly.
In fact, it is quite similar to that of the Schubert variety associated
with $w_0$ as in \cite{N} and then we obtain the following formula 
by the similar method in \cite{N}.
\begin{pro}
The action $e^c_i$ on 
$t\Theta^-_{\bfii0}(c_1,\cd,c_N)$ is given by
\[
e_i^c(t\Theta^-_{\bfii0}(c_1,\cd,c_N))
=t\Theta^-_{\bfii0}(c'_1,\cd,c'_N)
\]
where
\begin{equation}
c'_j\seteq 
c_j\cdot \frac{\displaystyle \sum_{1\leq m< j,\,i_m=i}
{c\cdot c_1^{a_{i_1,i}}\cd c_{m-1}^{a_{i_{m-1},i}}c_m}
+\sum_{j\leq m\leq N,\,i_m=i} 
{c_1^{a_{i_1,i}}\cd c_{m-1}^{a_{i_{m-1},i}}c_m}}
{\displaystyle \sum_{1\leq m\leq j,\,i_m=i} 
{c\cdot c_1^{a_{i_1,i}}\cd c_{m-1}^{a_{i_{m-1},i}}c_m}+
\mathop\sum_{j< m\leq N,\,i_m=i}  
{c_1^{a_{i_1,i}}\cd c_{m-1}^{a_{i_{m-1},i}}c_m}}.
\label{eici}
\end{equation}
The explicit forms of 
rational functions $\vep_i$ and $\gamma_i$ are:
\begin{equation}
 \vep_i(t\Theta^-_{\bfii0}(c))=
\left(\sum_{1\leq m\leq N,\,i_m=i} \frac{1}
{c_mc_{m+1}^{a_{i_{m+1},i}}\cd c_{N}^{a_{i_{N},i}}}\right)^{-1},\q
\gamma_i(t\Theta^-_{\bfii0}(c))
=\frac{\al_i(t)}{c_1^{a_{i_1,i}}\cd c_N^{a_{i_N,i}}}.
\label{th-vep-gamma}
\end{equation}
\end{pro}
{\sl Proof.}
If we rewrite $t\Theta^-_\bfii0(c)$ in the form,
\[
 t\cdot \al^\vee_{i_1}(c_1^{-1}\cd \al^\vee_{i_N}(c_N^{-1})
y_{i_1}(d_1)\cd y_{i_N}(d_N),
\]
then we easily get 
$d_m=c_mc_{m+1}^{a_{i_{m+1},i}}\cd c_{N}^{a_{i_{N},i}}$ for $m=1,\cd, N$.
Thus, we obtain the explicit form of $\vep_i$ as above.
To find the explicit form of the action $e_i^c$, the following formula
is crucial:
\[
x_i(a)y_j(b)=\begin{cases}
y_i(\frac{b}{1+ab})\al^\vee_i(1+ab)x_i(\frac{a}{1+ab})&\text{ if }i=j,\\
y_j(b)x_i(a)&\text{ if }i\ne j.
\end{cases}
\]
Applying this formula to (\ref{ei-action}) repeatedly, 
we have the above explicit 
action of $e_i^c$.\qed

\subsection{Ultra-Discretization of $\bbB_w=TB^-_w$}
Applying the ultra-discretization functor to
$\bbB_w$,
we obtain the free crystal ${\mathcal UD}(\bbB_{w_0})=X_*(T)\times \bbZ^{N}$,
where $N$ is the length of the longest element $w_0$.
Then define the map $\til h:{\mathcal UD}(\bbB_{w_0})
=X_*(T)\times \bbZ^{N}\,\,\to\,\,
X_*(T)(=P^*)$ as the projection to the left component and set
\[
 B_{w_0}(\lm^\vee):=\til h^{-1}(\lm^\vee),\qq
B_{f_B,\Theta^-_{\bfii0}}(\lm^\vee):=B_{w_0}(\lm^\vee)\cap 
B_{{f_B},\Theta^-_{\bfii0}},
\]
for $\lm^\vee\in X_*(T)=P^*$. Set $P^*_+:=\{h\in P^*|\Lm_i(h)\geq0
\text{ for any }i\in I\}$ and for $\lm^\vee=\sum_i\lm_ih_i$, we define
 $\lm=\sum_i\lm_i\Lm_i\in P_+$. Then, we have 
\begin{thm}[\cite{BK2}]
The set $B_{f_B,\Theta^-_{\bfii0}}(\lm^\vee)$ is non-empty if 
$\lm^\vee\in P^*_+$
 and in that case, $B_{f_B,\Theta^-_{\bfii0}}(\lm^\vee)$ is isomorphic to 
$B(\lm)^L$, which is the Langlands dual crystal associated with 
$\ge^L$. 
\end{thm}

It follows from 
(\ref{eici}) and (\ref{th-vep-gamma}) that we have 
\begin{thm}
\label{ud-tb}
Let $\lm^\vee\in P^*_+$. 
The explicit crystal structure of $B_{f_B,\Theta^-_\bfii0}(\lm^\vee)$
is as follows:
For $x=(x_1,\cd,x_N)\in B_{f_B,\Theta^-_\bfii0}(\lm^\vee)\subset \bbZ^N$, we
 have
\begin{equation}
\eit^n(x)=(x'_1,\cd,x'_N),
\end{equation}
where 
\begin{eqnarray}
&& x'_j=x_j+\min\left(\min_{1\leq m<j,i_m=i}(n+\sum_{k=1 }^m
 a_{i_k,i}x_k),
\min_{j\leq m\leq N,i_m=i}(\sum_{k=1 }^m a_{i_k,i}x_k)\right)
\label{eitn} \\
&&\qq -\min\left(\min_{1\leq m\leq j,i_m=i}(n+\sum_{k=1 }^m
 a_{i_k,i}x_k),
\min_{j<m\leq N,i_m=i}(\sum_{k=1 }^m a_{i_k,i}x_k)\right),
\nn
\\
&&\wt_i(x)=\lm(h_i)-\sum_{k=1}^Na_{i_k,i}x_k,
\label{wt-ud}\\
&&\vep_i(x)=\max_{1\leq m\leq
 N,i_m=i}(x_m+\sum_{k=m+1}^{N}a_{i_k,i}x_k),
\label{vep-ud}
\end{eqnarray}
and $x=(x_1,\cd,x_N)$ belongs to $B_{f_B,\Theta^-_\bfii0}(\lm^\vee)$ if and
 only if ${\mathcal UD}(f_B)(x)\geq0$.
\end{thm}
It follows immediately from (\ref{eitn}):
\begin{lem}
\label{lem-ef}
Set $X_m:=\sum_{k=1 }^m a_{i_k,i}x_k$, ${\mathcal X}^{(i)}:=\min\{X_m|1\leq
 m\leq N,i_m=i\}$ $(i\in I)$ and $M^{(i)}:=\{l|1\leq l\leq N,i_l=i,
 X_l={\mathcal X}^{(i)}\}$. 
Define $m_e:=\max(M^{(i)})$ and $m_f:=\min(M^{(i)})$:
for $x\in B_{f_B,\Theta^-_\bfii0}(\lm^\vee)$, we get 
\begin{eqnarray}
&&\eit(x)=\begin{cases}(x_1,\cd,x_{m_e}-1,\cd,x_N)
&\text{ if }{\mathcal UD}(f_B)(x_1,\cd,x_{m_e}-1,\cd,x_N)\geq0,\\
0&\text{otherwise,}
\end{cases}
\label{th-eaction}\\
&&\fit(x)=\begin{cases}(x_1,\cd,x_{m_f}+1,\cd,x_N)
&\text{ if }{\mathcal UD}(f_B)(x_1,\cd,x_{m_f}+1,\cd,x_N)\geq0,\\
0&\text{otherwise.}
\end{cases}
\label{th-faction}
\end{eqnarray}
\end{lem}
Finally, due to the results in Sect.2 and in this section, 
we obtain the following theorem 
\begin{thm}
\label{coin-thm}
If we have $B_{f_B,\Theta^-_\bfii0}(\lm^\vee)
=\Sigma_{\bfii0^{-1}}[\lm]^L$ as a set.
Then they are isomorphic each other as crystals, where ${}^L$ means the
Langlands dual crystal, that is, it is defined by the transposed Cartan 
matrix and $\bfii0^{-1}$ means the opposite order of $\bfii0$.
\end{thm}
{\sl Proof.}
The coincidence of the actions $\eit$ and $\fit$ are shown by comparing
(\ref{action-f}) and  (\ref{action-e}) with  (\ref{th-eaction}) and 
(\ref{th-faction}) since the following are equivalent:\\
(a) $X_k$ is the minimum.\\
(b) $\sigma^{(i)}(x)=\sigma^{(0)}_i(x)+\lm_i-X_k$ is the maximum.\\
Similarly, comparing (\ref{wt-vep-vp-1}) with (\ref{wt-ud}) and 
(\ref{wt-vep-vp4}) with (\ref{vep-ud}) respectively, 
we obtain the coincidence 
of $\vep_i$ and $wt_i$. \qed
\renewcommand{\thesection}{\arabic{section}}
\section{Explicit form of the decoration $f_B$ of type $A_n$}
\setcounter{equation}{0}
\renewcommand{\theequation}{\thesection.\arabic{equation}}

\subsection{Generalized Minors and the function $f_B$}

For this subsection, see \cite{BFZ,BZ,BZ2}.
Let $G$ be a simply connected simple algebraic groups over $\bbC$ and 
$T\subset G$ a maximal torus. 
Let  $X^*(T):=\Hom(T,\bbC^\times)$ and $X_*(T):=\Hom(\bbC^\times,T)$ be
the lattice of characters and co-characters respectively.
We identify $P$ (resp. $P^*$) with $X^*(T)$ 
(resp. $X_*(T)$) as above.
\begin{df}
For $\mu\in P_+$, the
{\it principal minor} $\Del_\mu:G\to\bbC$ is defined as
\[
 \Del_\mu(u^-tu^+):=\mu(t)\q(u^\pm\in U^\pm,\,\,t\in T).
\]
Let $\gamma,\del\in P$ be extremal weights such that 
$\gamma=u\mu$ and $\del=v\mu$ for some $u,v\in W$. 
Then the {\it generalized minor} $\Del_{\gamma,\del}$ is defined
by
\[
 \Del_{\gamma,\del}(g):=\Del_\mu(\ovl u^{-1}g\ovl v)
\q(g\in G),
\]
which is a regular function on $G$.
\end{df}
\begin{lem}[\cite{BK2}]
Suppose that $G$ is simply connected.
\begin{enumerate}
\item
For $u\in U$ and $i\in I$, we have 
$\Del_{\mu,\mu}(u)=1$ and $\chi_i(u)=\Del_{\Lm_i,s_i\Lm_i}(u)$,
where $\Lm_i$ be the $i$th fundamental weight.
\item
Define the map $\pi^+:B^-\cdot U\to U$ by $\pi^+(bu)=u$ for 
$b\in B^-$ and $u\in U$. For any $g\in G$, we have 
\begin{equation}
 \chi_i(\pi^+(g))=\frac{\Del_{\Lm_i,s_i\Lm_i}(g)}
{\Del_{\Lm_i,\Lm_i}(g)}.
\end{equation}
\end{enumerate}
\end{lem}
\begin{pro}[\cite{BK2}]
The function $f_B$ in (\ref{f_B}) is described as follows:
\begin{equation}
 f_B(g)=\sum_i\frac{\Del_{w_0\Lm_i,s_i\Lm_i}(g)
+\Del_{w_0s_i\Lm_i,\Lm_i}(g)}
{\Del_{w_0\Lm_i,\Lm_i}(g)}
\end{equation}
\end{pro}
Let ${\bf i}=i_1\cd i_N$ be a reduced word for the longest Weyl
group element $w_0$. 
For $t\Theta_{\bf i}^-(c)\in \bbB_{w_0}=T\cdot B^-_{w_0}$, 
we get the following formula.
\begin{equation}
 f_B(t\Theta_{\bf i}^-(c))
=\sum_i\Del_{w_0\Lm_i,s_i\Lm_i}(\Theta_{\bf i}^-(c))
+\al_i(t)\Del_{w_0s_i\Lm_i,\Lm_i}(\Theta_{\bf i}^-(c)).
\label{fb-th}
\end{equation}
\subsection{Bilinear Forms}

Let $\omega:\ge\to\ge$ be the anti involution 
\[
\omega(e_i)=f_i,\q
\omega(f_i)=e_i\,\q\omega(h)=h,
\] and extend it to $G$ by setting
$\omega(x_i(c))=y_i(c)$, $\omega(y_i(c))=x_i(c)$ and $\omega(t)=t$
$(t\in T)$.

There exists a $\ge$(or $G$)-invariant bilinear form on the
finite-dimensional  irreducible
$\ge$-module $V(\lm)$ such that 
\[
 \lan au,v\ran=\lan u,\omega(a)v\ran,
\q\q(u,v\in V(\lm),\,\, a\in \ge(\text{or }G)).
\]
For $g\in G$, 
we have the following simple fact:
\[
 \Del_{\Lm_i}(g)=\lan gu_{\Lm_i},u_{\Lm_i}\ran.
\]
Hence, for $w,w'\in W$ we have
\begin{equation}
 \Del_{w\Lm_i,w'\Lm_i}(g)=
\Del_{\Lm_i}({\ovl w}^{-1}g\ovl w')=
\lan {\ovl w}^{-1}g\ovl w'\cdot u_{\Lm_i},u_{\Lm_i}\ran
=\lan g\ovl w'\cdot u_{\Lm_i}\, ,\, \ovl{w}\cdot u_{\Lm_i}\ran,
\label{minor-bilin}
\end{equation}
where $u_{\Lm_i}$ is a properly normalized highest weight vector in
$V(\Lm_i)$ and note that $\omega(\ovl s_i^{\pm})=\ovl s_i^{\mp}$.

\subsection{Explicit form of  $f_B(t\Theta_{\bf i}^-(c))$ of type $A_n$}
\label{fb-An}

Now, we consider the type $A_n$, that is, $G=SL_{n+1}(\bbC)$. 
We fix the reduced longest word ${\bf  i_0}=\underbrace{1,2,\cd,n}_{},
\underbrace{1,2,\cd,n-1}_{},\cd,\underbrace{1,2,3}_{},1,2,1$.
This is just the opposite order $\io_0=\bfii0^{-1}$ as in Sect.2.
To obtain the explicit form of $f_B(t\Theta_{\bf i_0}^-(c))$, 
by (\ref{fb-th}) it suffices to know 
$\Del_{w_0\Lm_j,s_j\Lm_j}(\Theta_{\bf j_0}^-(c))$ and 
$\Del_{w_0s_j\Lm_j,\Lm_j}(\Theta_{\bf j_0}^-(c))$ for 
\[
c=(c_{i,j}|i+j\leq n+1)=
(c_{1,1},c_{1,2},\cd,c_{1,n},c_{2,1},c_{2,2},\cd,c_{2,n-1},\cd
c_{n-1,1},c_{n-1,2},c_{n,1})\in (\bbC^\times)^N.
\]

\begin{thm}
\label{thm-a}
For $c\in (\bbC^\times)^N$ as above, 
we have the following explicit forms:
\begin{eqnarray}
&&\Del_{w_0\Lm_j,s_j\Lm_j}(\Theta_\bfii0^-(c))=
c_{n-j+1,1}+\frac{c_{n-j+1,2}}{c_{n-j+2,1}}+\frac{c_{n-j+1,3}}{c_{n-j+2,2}}
+\cd+\frac{c_{n-j+1,j}}{c_{n-j+2,j-1}}, 
\label{del1}\\
&&\Del_{w_0s_j\Lm_j,\Lm_j}(\Theta_\bfii0^-(c))=
\frac{1}{c_{j,1}}+
\frac{c_{j-1,1}}{c_{j-1,2}}+\frac{c_{j-2,2}}{c_{j-2,3}}
+\cd+\frac{c_{1,j-1}}{c_{1,j}},\,\,
(j\in I).
\label{del2}
\end{eqnarray}
\end{thm}
The proof of this theorem will be given in the next subsection.
\subsection{Proof of Theorem \ref{thm-a}}

Let $V_1:=V(\Lm_1)$ be the vector representation of 
$\ssl_{n+1}(\bbC)$ with
the standard basis $\{v_1,\cd,v_{n+1}\}$, and 
$\{e_i,f_i,h_i\}_{i=1,\cd,n}$ the 
Chevalley generators of $\ssl_{n+1}(\bbC)$.
Their actions on the basis vectors are as follows:
\begin{equation}
e_iv_j=\begin{cases}v_{i}&\text{ if }j=i+1,\\
0&\text{otherwise},
\end{cases}
\q
f_iv_j=\begin{cases}v_{i+1}&\text{ if }j=i,\\
0&\text{otherwise},
\end{cases}
\q
h_iv_j=\begin{cases}v_{i}&\text{ if }j=i,\\
-v_{i-1}&\text{ if }j=i-1\text { and }i\ne1,\\
0&\text{otherwise},
\end{cases}
\end{equation}
By these explicit actions 
we know that $e_i^2=f_i^2=0$ on $V_1$.
Thus, we can write 
${\pmb x}_i(c):=\al_i^\vee(c^{-1})x_i(c)=c^{-h_i}(1+c\cdot e_i)$ and 
${\pmb y}_i(c):=y_i(c)\al_i^\vee(c^{-1})=(1+c\cdot f_i)c^{-h_i}$ on $V_1$
and then 
\begin{equation}
\pmbx_i(c)v_j=\begin{cases}cv_{i+1}+v_i&\text{ if }j=i+1,\\
cv_i&\text{ if }j=i,\\
v_j&\text{ otherwise, }
\end{cases}\q
\pmby_i(c)v_j=\begin{cases}c^{-1}v_i+v_{i+1}&\text{ if }j=i,\\
cv_i&\text{ if }j=i-1,\\
v_j&\text{ otherwise. }
\end{cases}
\end{equation}

For $c=(c_1,\cd,c_i)\in (\bbC^\times)^i$ set
$X^{(i)}(c):=\pmbx_i(c_i)\cd \pmbx_1(c_1)$ and its action on the basis vector is
as follows:
\begin{equation}
X^{(i)}(c)v_k=\begin{cases}
c_k^{-1}c_{k-1}v_k+v_{k-1}&\text{ if }k<i+1,\\
c_iv_{i+1}+v_i&\text{ if }k=i+1,\\
v_k&\text{ if }k>i+1.
\end{cases}
\end{equation}
Now, for $c=(c_{k,i})_{1\leq i,k\leq n,i+k\leq n+1}
\in(\bbC^\times)^{\frac{n(n+1)}{2}}$ set 
$c^{(k)}:=(c_{n+1-k,k},c_{n+1-k,k-1},
\cd,c_{n+1-k,2},c_{n+1-k,1})$ and 
\[
 X(c):=X^{(1)}(c^{(1)})X^{(2)}(c^{(2)})\cd X^{(n-1)}(c^{(n-1)})X^{(n)}(c^{(n)})
\]
Here, note that 
\begin{equation}
\omega(\Theta^-_{{\bf i}_0}(c))=X(c).
\label{omega-th-x}
\end{equation}
Writing 
\[
X(c)v_i=\sum_{k=1}^i\xx(i,k)v_k,
\]
we shall get the explicit form of the coefficient $\xx(i,k)$ 
with the direct calculations:
For $i=1,\cd,n$ and $k=1,\cd,i$, set 
\[
 {}^i{\bf m}_k:=\{M|M\subset\{1,\cd,n-k+1\},\sharp M=n-i+1\}.
\]
For $M\in {}^i{\bf m}_k$, write $M=M_1\sqcup \cd\sqcup M_{i-k+1}$ where 
each $M_j$ ($j=1,\cd,s:=i-k+1$) is a consecutive subsequence of $M$ satisfying
$\min(M_{l})-\max(M_j)=l-j+1$ for any $1\leq j<l\leq s$ 
if  both $M_j$ and $M_l$ are non-empty,
which is called a segment of $M$. 
For $M=M_1\sqcup\cd M_s\in {}^i{\bf m}_k$, write each segment:
\[
M_1=\{1,2,\cd,j_1-1\},\,\
M_2=\{j_1+1,j_1+2,\cd,j_2-1\},\,\,
M_s=\{j_{i-k}+1,j_{i-k}+2,\cd,n-k+1\},
\]
where $1\leq j_1<j_2<\cd<j_{i-k}\leq n-k+1$ and set 
\[
 c^M:=\frac{c_{1,i-1}\cd c_{j_1-1, i-1}\cd c_{j_{i-k}+1,k-1}
\cd c_{n-k+1,k-1}}
{c_{1,i}\cd c_{j_1-1,i}\cd c_{j_{i-k}+1,k}
\cd c_{n-k+1,k}}.
\]
\begin{pro}
We have the following explicit form of $\xx(i,k)$.
\begin{equation}
\xx(i,k)=\sum_{M\in{}^i{\bf m}_k}c^M.
\end{equation}
\end{pro}
This formula is obtained by direct calculations.

For the module $V(\Lm_j)$ ($j>1$), let us denote 
its normalized highest (resp. lowest)weight vector by $u_{\Lm_j}$ 
(resp. $v_{\Lm_j}$). 
Set 
\begin{eqnarray*}
&&[i_1,\cd,i_j]:=v_{i_1}\wedge v_{i_2}\wedge\cd\wedge v_{i_j}
\in \bigwedge^jV_1 \\
&&I_j:=\{[i_1,i_2,\cd,i_j]\,|\,1\leq i_1<i_2<\cd<i_j\leq n+1\}.
\end{eqnarray*}
$I_j$ is a normal basis of $V(\Lm_j)$ with the weight 
$\sum_{k=1}^j(\Lm_{i_k}-\Lm_{i_k-1})$.
Indeed, $u_{\Lm_j}=v_1\wedge v_2\wedge \cd\wedge v_j$ and 
$v_{\Lm_j}=v_{n-k}\wedge v_{n-k+1}\wedge \cd\wedge v_{n+1}$.
The actions of $e_i$ and $f_i$ on the vector $[i_1,\cd,i_j]$ are given by
\begin{eqnarray}
&&e_i[i_1,\cd,i_j]=\begin{cases}
[i_1,\cd,i_{k-1},i,i_{k+1},\cd,i_j]&\text{ if }i_k=i+1, \,\,
i_{k-1}<i \text{ for some }k,\\
0&\text{otherwise},
\end{cases}\\
&&f_i[i_1,\cd,i_j]=\begin{cases}
[i_1,\cd,i_{k-1},i+1,i_{k+1},\cd,i_j]&\text{ if }i_k=i, \,\,
i_{k+1}>i+1 \text{ for some }k,\\
0&\text{otherwise}.
\end{cases}
\end{eqnarray}
It follows from the formula (\ref{minor-bilin}) and (\ref{omega-th-x}) that 
\begin{equation}
\Del_{w_0\Lm_j,s_i\Lm_j}(\Theta_{\bf i_0}^-(c))
=\lan \Theta_{\bf i_0}^-(c)\ovl s_j\cdot u_{\Lm_j}, \ovl w_0\cdot
u_{\Lm_j}\ran
=\lan\ovl s_j\cdot u_{\Lm_j}, X(c)\ovl w_0\cdot u_{\Lm_j}\ran.
\end{equation}
Since $\ovl s_j\cdot u_{\Lm_j}=[1,2,\cd,j-1,j+1]$ and 
$\ovl w_0\cdot u_{\Lm_j}=v_{\Lm_j}$, to obtain 
$\Del_{w_0\Lm_j,s_i\Lm_j}(\Theta_{\bf i_0}^-(c))$ it suffices to find
the coefficient of $[1,2,\cd,j-1,j+1]$ in $X(c)v_{\Lm_j}$.
\begin{lem}
We have  
\begin{equation}
X^{(j+1)}(c^{(j+1)})\cd X^{(n)}(c^{(n)})v_{\Lm_j}
=[2,3,\cd,j,j+1]+\sum_{1\leq i_1<\cd<i_j>j+1}c_{i_1,\cd,i_j}[i_1,\cd,i_j],
\end{equation}
where $c_{i_1,\cd,i_j}\in\bbC$ is the coefficient.
\end{lem}
{\sl Proof.}
First, let us see 
$W_n:=X^{(n)}(c^{(n)})v_{\Lm_j}$. It is easily to see that 
\begin{eqnarray}
&&W_n=\pmbx_n(c_{1,n})\cd \pmbx_1(c_{1,1})v_{\Lm_j}
\label{Xn} \\
\qq\qq &&=
[n+1-j,n+2-j,\cd,n-1,n]+
\sum_{1\leq i_1<\cd<i_{j-1}<n+1}c_{i_1,\cd,i_{j-1},n+1}[i_1,\cd,i_{j-1},n+1],
\nn
\end{eqnarray}
since the term 
$c_{1,n}^{-h_n}c_{1,n}e_n\cdot
c_{1,n-1}^{-h_{n-1}}c_{1,n-1}e_{n-1}
\cd c_{1,n+1-j}^{-h_{n+1-j}}c_{1,n+1-j}e_{n+1-j}\cdot
c_{1,n-j}^{-h_{n-j}}\cd
c_{1,1}^{-h_1}$ in $X^{(n)}(c^{(n)})$
gives the leading term $[n+1-j,n+2-j,\cd,n-1,n]$ in (\ref{Xn}).
Indeed, the basis vectors appearing in
$X^{(n-1)}(c^{(n-1)})[i_1,\cd,n+1]$
are in  the form $[\cd\cd,n+1]$, which means that 
we may see only the vector
$X^{(n-1)}(c^{(n-1)})[n+1-j,n+2-j,,\cd,n-1,n]$ in 
$X^{(n-1)}(c^{(n-1)})X^{(n)}(c^{(n)})v_{\Lm_j}$.
By considering similarly, we obtain 
\[
 X^{(n-1)}(c^{(n-1)})X^{(n)}(c^{(n)})v_{\Lm_j}
=[n-j,n+1-j,\cd,n-2,n-1]+
\sum_{1\leq i_1<\cd<i_{j}\geq n}c_{i_1,\cd,i_j}[i_1,\cd,i_j].
\]
Repeating this process, we get the desired result.\qed

We shall see the action  of $\ovl X_j:=
X^{(1)}(c^{(1)})\cd X^{(j)}(c^{(j)})$
on the vector $[2,3,\cd,j,j+1]$.
The following lemma is shown easily.
\begin{lem}
In the expansion of 
\[
 \ovl X_j=c_{n,1}^{-h_n}(1+c_{n,1}e_1)\cd c_{j,1}^{-h_1}
(1+c_{j,1}e_1),
\]
the only terms $E_m:=A\cdot B_m\cdot C_m\cdot D_m$ ($m=0,\cd,j-1$)
produce the vector $[1,2,\cd,j-1,j+1]$ in 
$X^{(1)}(c^{(1)})\cd X^{(j)}(c^{(j)})[2,3,\cd,j,j+1]$, where 
\begin{eqnarray*}
&&A:=\ch(1,n)\ch(2,n-1)\ch(1,n-1)\cd 
\ch(j-2,n-j+3)\ch(j-3,n-j+3)\cd \ch(1,n-j+3),\\
&&B_m:=\ch(j-1,n-j+2)\cc(j-1,n-j+2)e_{j-1}\cd 
\ch(m+2,n-j+2)\cc(m+2,n-j+2)e_{m+2}\cdot 
\ch(m+1,n-j+2)\cc(m+1,n-j+2)e_{m+1},\\
&&C_m:=\ch(m,n-j+2)\cd\ch(2,n-j+2)\ch(1,n-j+2)\ch(j,n-j+1)\ch(j-1,n-j+1)
\cd\ch(m+1,n-j+1),\\
&&D_m:=\ch(m,n-j+1)\cc(m,n-j+1)e_{m}\cdots
\ch(2,n-j+1)\cc(2,n-j+1)e_{2}\cdot 
\ch(1,n-j+1)\cc(1,n-j+1)e_{1},
\end{eqnarray*}
where we understand $D_0=1$.
\end{lem}
For $m=1,2,\cd,j-1$  it is trivial that 
\[
C_m\cdot D_m[2,3,\cd,j,j+1]=
\frac{\cc(m+1,n-j+1)}{\cc(m,n-j+2)}
[1,2,\cd,m,m+2,m+3,\cd,j,j+1].
\]
And then 
\[
 A\cdot 
B_m\cdot C_m\cdot D_m[2,3,\cd,j,j+1]=
\frac{\cc(m+1,n-j+1)}{\cc(m,n-j+2)}
[1,2,3,\cd,j-1,j+1].
\]
For $m=0$, we have 
\[
 A\cdot B_0\cdot C_0\cdot D_0[2,3,\cd,j,j+1]=\cc(1,n-j+1)[1,2,\cd,j-1,j+1].
\]
Finally, we obtain that 
the coefficient of $[1,2,\cd,j-1,j+1]$ in $X(c)v_{\Lm_j}$ is 
\begin{equation}
\cc(1,n-j+1)+\sum_{m=1}^{j-1}\frac{\cc(m+1,n-j+1)}{\cc(m,n-j+2)},
\end{equation}
which is just $\Del_{w_0\Lm_j,s_i\Lm_j}(\Theta_{\bf i_0}^-(c))$ 
and then we have shown (\ref{del1}) in Theorem \ref{thm-a}. 
The formula (\ref{del2}) would be shown by the similar way to 
(\ref{del1}). \qed

Note that for $j=1$ and $k=2$, we find 
$\xx(1,2)=\Del_{w_0\Lm_1,s_1\Lm_1}(\Theta^-_{{\bf i}_0}(c))$.

\renewcommand{\thesection}{\arabic{section}}
\section{Ultra-Discretization and Polyhedral Realizations of type $A_n$}
\setcounter{equation}{0}
\renewcommand{\theequation}{\thesection.\arabic{equation}}



In this section, we shall only treat the type $A_n$. Then we identify 
$P$ with $P^*$ by $\lm\leftrightarrow \lm^\vee$.
Let us describe the explicit form of $B_{f_B,\Theta^-_{\bfii0}}(\lm)$ for
type $A_n$ applying the result in Theorem \ref{ud-tb} and 
show the coincidence of the crystals $B_{f_B,\Theta^-_{\bfii0}}(\lm)$
and 
$\Sigma_{\io_0}[\lm]$ in Sect.2 using Theorem \ref{coin-thm}.

For $\ge=\ssl_{n+1}(\bbC)$, 
let $\bfii0$ be as in \ref{fb-An}. Then 
we have the following:
\begin{lem}
\label{lem-coin-An}
The crystal $B_{f_B,\Theta^-_{\bfii0}}(\lm)$ is defined 
by 
\begin{equation}
B_{f_B,\Theta^-_{\bfii0}}(\lm):=
\left\{(x_{k,l}|k+l\leq n+1)\in \bbZ^{N}\,\,
\begin{array}{|l}
\hbox{
$x_{1,i} \geq x_{2,i-1} \geq \cd \geq x_{i,1} \geq 0$
for $1 \leq i \leq n$}
\\
\hbox{
$\lm_i\geq x_{j,i-j+1}-x_{j,i-j}$ for
$1\leq j\leq i\leq n$.}
\end{array}\right\},
\label{explicit-bfb}
\end{equation}
where $N=\frac{n(n+1)}{2}$.
\end{lem}
{\sl Proof.}
We shall see the explicit form of ${\mathcal UD}(f_B)(x)$. Indeed, 
by virtue of (\ref{fb-th}), it is sufficient to know 
the forms of $\Del_{w_0\Lm_j,s_j\Lm_j}(\Theta_\bfii0^-(c))$ 
and $\Del_{w_0s_j\Lm_j,\Lm_j}(\Theta_\bfii0^-(c))$, which are given in 
(\ref{del1}) and (\ref{del2}). Thus, 
we have
\[
 {\mathcal UD}(f_B)(t,x)=\min_{1\leq j\leq n}
({\mathcal UD}(\Del_{w_0\Lm_j,s_j\Lm_j}(\Theta_\bfii0^-))(x),
{\mathcal UD}(\al_j(t))
+{\mathcal UD}(\Del_{w_0s_j\Lm_j,\Lm_j}(\Theta_\bfii0^-))(x))
\]
and 
\begin{eqnarray}
&&{\mathcal UD}(\Del_{w_0\Lm_j,s_j\Lm_j}(\Theta_\bfii0^-))(x)
=\min_{k=1,\cd,j}(x_{n-j+1,k}-x_{n-j+2,k-1}),\\
&&{\mathcal UD}(\Del_{w_0s_j\Lm_j,\Lm_j}(\Theta_\bfii0^-))(x))
=\min_{k=1,\cd,j}(x_{j-k+1,k-1}-x_{j-k+1,k}),
\end{eqnarray}
where $x_{j,k}={\mathcal UD}(c_{j,k})$ and 
we understand $x_{m,0}=0$.
Hence, if we identify ${\mathcal UD}(\al_j(t))$ with $\lm_j$, then 
the condition ${\mathcal UD}(f_B)(\lm,x)\geq0$ 
in Theorem \ref{ud-tb}
is equivalent to the condition in 
(\ref{explicit-bfb}).\qed

\begin{thm}
\label{thm-b}
For any dominant integral weight $\lm$, 
there exists the following isomorphism of crystals 
 $B_{f_B,\Theta^-_{\bfii0}}[\lm]\cong\Sigma_{\io_0}[\lm]$
where $\Sigma_{\io_0}[\lm]$ is as in Corollary \ref{cor-A_n} and 
$\io_0=\bfii0^{-1}$.
\end{thm}
{\sl Proof.} By Theorem \ref{coin-thm} it is necessarily for us to show
that $B_{f_B,\Theta^-_{\bfii0}}[\lm]=\Sigma_{\io_0}[\lm]$ as a set, 
which is shown by Corollary \ref{cor-A_n} and Lemma \ref{lem-coin-An}.
\qed

\renewcommand{\thesection}{\arabic{section}}
\section{Elementary Characters and Monomial Realization of Crystals}
\setcounter{equation}{0}
\renewcommand{\theequation}{\thesection.\arabic{equation}}

We shall see the elementary characters as in Sect4 from the different point of
view, that is, the monomial realization of crystals.

Let us introduce the monomial realization of crystals (See \cite{K7,Nj}).
For variables $\{Y_{m,i}|i\in I, m\in\bbZ.\}$, define the set of monomials
\[
 \cY:=\{Y=\prod_{m\in\bbZ,i\in I}
Y_{m,i}^{l_{m,i}}|l_{m,i}\in \bbZ\setminus\{0\}\text{
 except for finitely many }(m,i)\}.
\]
Fix a set of integers $p=(p_{i,j})_{i,j\in I,i\ne j}$ such that 
$p_{i,j}+p_{j,i}=1$. For this $p:=(p_{i,j})_{i,j\in I,i\ne j}$ and a
generalized Cartan matrix $(a_{i,j})_{i,j\in I}$, set
\[
 A_{m,i}=Y_{m,i}Y_{m+1,i}\prod_{j\ne i}Y_{m+p_{j,i},j}^{a_{j,i}}.
\]
Note that 
for any cyclic order $\io=\cd (i_1i_2\cd i_n)(i_1i_2\cd i_n)\cd$
s.t. $\{i_1,\cd,i_n\}=I$,
we can associate the following $(p_{i,j})$ by:
\[
 p_{i_a,i_b}=\begin{cases}
1&a<b,\\0&a>b.\end{cases}
\]
For example, if we take $\io=\cd (213)(213)\cd$, then we have
$p_{2,1}=p_{1,3}=p_{2,3}=1$ and $p_{1,2}=p_{3,1}=p_{3,2}=0$.
Thus, we can identify a cyclic order $\cd(i_1\cd i_n)(i_1\cd i_n)\cd$ 
with such $(p_{i,j})$.

For a monomial $Y=\prod_{m,i}Y_{m,i}^{l_{m,i}}$, 
set 
\begin{eqnarray*}
&&\hspace{-20pt}wt(Y)=\sum_{i,m} l_{m,i}\Lm_i, \,\,
\vp_i(Y)=\operatorname{max}_{k\in\bbZ}\{\sum_{k\leq m}l_{m,i}\},\,\,
\vep_i(Y)=\vp_i(Y)-wt(Y)(h_i),\\
&&\hspace{-20pt}\fit(Y)=\begin{cases}
A_{n_f,i}^{-1}\cdot Y&\text{ if }\vp_i(Y)>0,\\
0&\text{ if }\vp_i(Y)=0,
\end{cases}\q\q
\eit(Y)=\begin{cases}
A_{n_e,i}\cdot Y&\text{ if }\vep_i(Y)>0,\\
0&\text{ if }\vep_i(Y)=0,
\end{cases}\\
&&n_f=\min\{n|\vp_i(Y)=\sum_{k\leq n}m_{k,i}\},\q
n_e=\max\{n|\vp_i(Y)=\sum_{k\leq n}m_{k,i}\}.
\end{eqnarray*}

\begin{thm}[\cite{K7,Nj}]
\begin{enumerate}
\item
In the above setting, $\cY$ is a crystal, which is denoted by $\cY(p)$.
\item
If $Y\in\cY(p)$ satisfies $\vep_i(Y)=0$ for any $i\in I$, 
then the connected component containing $Y$ is isomorphic to
$B(wt(Y))$.
\end{enumerate}
\end{thm}

In the above setting, for type $A_n$ take $(p_{i,j})_{i,j\in I,i\ne j}$
such that $p_{i,j}=1$ for $i<j$, $p_{i,j}=0$
for $i>j$, which corresponds to the cyclic order
${\mathbf i}=(12\cd n)(12\cd n)\cd$.
Then we obtain 
\begin{pro}
\label{mono-cry}
The crystal containing the monomial $Y_{n-i+1,1}$ (resp. $Y_{i,1}^{-1}$)
is isomorphic to $B(\Lm_1)$ (resp. $B(\Lm_n)$)
and all basis vectors are given by
\begin{eqnarray*}
&&
\til f_k\cd \til f_2\til
f_1(Y_{n-i+1,1})=\frac{Y_{n-i+1,k+1}}{Y_{n-i+2,k}}\in B(\Lm_1),\\
&&\til e_k\cd\til e_2\til e_1(Y_{i,1}^{-1})
=\frac{Y_{i-k,k}}{Y_{i-k,k+1}}\in B(\Lm_n)\q
(k=1,\cd,n).
\end{eqnarray*}
\end{pro}
{\sl Proof.}
The explicit form of $A_{m,i}$ $(m\in\bbZ,i\in I)$ is as follows:
\begin{equation}
A_{m,i}=\begin{cases}
Y_{m,1}Y_{m,2}^{-1}Y_{m+1,1}&\text{ if }i=1,\\
Y_{m,i}Y_{m,i+1}^{-1}Y_{m+1,i-1}^{-1}Y_{m+1,i}&\text{ if }i\ne1,n,\\
Y_{m,n}Y_{m+1,n-1}^{-1}Y_{m+1,n}&\text{ if }i=n.
\end{cases}
\end{equation}
Then, applying $\eit$ and $\fit$ repeatedly, we obtain the results.
For example, 
\[
\til f_1(Y_{n-i+1,1})=Y_{n-i+1,1}\cdot A_{n-i+1,1}^{-1}
=\frac{Y_{n-i+1,2}}{Y_{n-i+2,1}}.
\]
\qed

Applying this results to Theorem \ref{thm-a} and 
changing the variable $Y_{m,l}$ to $c_{m,l}$, we find:
\begin{pro}
\label{del-mono}
For $j=1,\cd,n$  we have
\begin{eqnarray*}
&&\chi_j(\pi^+(w_0^{-1}t\Theta_\bfi^-(c)))=
\Del_{w_0\Lm_j,s_j\Lm_j}(\Theta_\bfi^-(c))=
\sum_{k=0}^{j-1}\til f_k\cd \til f_2\til f_1(c_{n-j+1,1}),\\
&&\chi_j(\pi^+(w_0^{-1}\eta(t\Theta_\bfi^-(c))))=
\al_j(t)\Del_{w_0s_j\Lm_j,\Lm_j}(\Theta_\bfi^-(c))=
\al_j(t)
\sum_{k=0}^{j-1}\til e_k\cd\til e_2\til e_1(c_{j,1}^{-1}).
\end{eqnarray*}
\end{pro}
Note that 
$\{\til f_k\cd \til f_2\til f_1(c_{n-i+1,1})|0\leq k<i\}
=B(\Lm_1){\tiny s_{k-1}\cd s_2s_1}$
is the Demazure crystal associated with the Weyl group element 
$s_{k-1}\cd s_2s_1$ (\cite{K3}).

Observing Proposition \ref{del-mono}, we present the following
conjecture:
\begin{conj}
\label{conj1}
There exists certain reduced longest word 
$\bfi=(i_1,\cd, i_N)$ and $p=(p_{i,j})_{i\ne j}$ 
such that for any $i\in I$, 
there exist Demazure crystal $B^-_{w}(i)\subset B(\Lm_k)$,
Demazure crystal $B^+_{w'}(i)\subset B(\Lm_j)$ and 
positive integers $\{a_b,a_{b'}|b\in B^-_w,b'\in B^+_{w'}\}$ satisfying
\begin{eqnarray*}
&&\chi_i(\pi^+(w_0^{-1}t\Theta_\bfi^-(c)))=
\Del_{w_0\Lm_i,s_i\Lm_i}(\Theta_\bfi^-(c))=
\sum_{b\in B^-_w(i)}a_b
m_b(c),\\
&&\chi_i(\pi^+(w_0^{-1}\eta(t\Theta_\bfi^-(c))))=
\al_i(t)\Del_{w_0s_i\Lm_i,\Lm_i}(\Theta_\bfi^-(c))=
\al_j(t)\sum_{b'\in B^+_{w'}(i)}a_{b'}m_{b'}(c),
\end{eqnarray*}
where $m_b(c)\in\cY(p)$ is the monomial corresponding to $b\in B(\Lm_k)$
associated with $p=(p_{i,j})_{i\ne j}$.
\end{conj}

We would see the answers to this
 conjecture for other type of Lie algebras in the 
forthcoming papers.

Suppose that this conjecture is right and then we can deduce the
following: 
\begin{cor}
In the setting of the above conjecture, define the linear function $\til
m_b(x):={\mathcal UD}(m_b)(x)$ ($x\in \bbZ^N$) and set
\[
 \wtil\Sigma_{{\bf i}^{-1}}[\lm]:=
\{x=(x_N,\cd,x_1)\in \bbZ^N_{{\bf i}^{-1}}[\lm]\,|\,
\til m_b\geq0,\,\,\lm_i+\til m_{b'}\geq0\text{ for any }b\in
B^-_w(i),b'\in B^+_{w'}(i)\,\,(i\in I)
\}.
\]
Then this is equipped with the crystal structure and 
isomorphic to the crystal $B(\lm)$.
\end{cor}
{\sl Proof.} Since ${\mathcal UD}(f_B)(\lm,x)\geq0$ is equivalent to the 
condition of the set $ \wtil\Sigma_{{\bf i}^{-1}}[\lm]$, we know that 
the set $ \wtil\Sigma_{{\bf i}^{-1}}[\lm]$ coincides with 
$B_{f_B,\Theta^-_\bfii0}(\lm^\vee)$.\qed

We call $\wtil\Sigma_{{\bf i}^{-1}}[\lm]$ the {\it refined polyhedral
realization} associated with the monomial realizations $\cY(p)$.

\bibliographystyle{amsalpha}

\end{document}